\documentclass[10pt,a4paper]{article}
\usepackage{amsfonts}
\usepackage{amssymb}
\usepackage{amsmath}
\usepackage{epsfig}
\usepackage{amsthm}
\usepackage[latin1]{inputenc}
\usepackage[english]{babel}

\newtheorem{defi}{Definition}[section]
\newtheorem{proposition}[defi]{Proposition}
\newtheorem{theoreme}[defi]{Theorem}

\newtheorem{lemme}[defi]{Lemma}
\newtheorem{corollaire}[defi]{Corollary}
\newtheorem{remarque}[defi]{Remark}

\newcommand{\F}{\mathbb F} 
\newcommand{\Z}{\mathbb Z}
\newcommand{\N}{\mathbb N}

\newcommand{\C}{\mathbb C}
\newcommand{\p}{\mathbb P}
\newcommand{\A}{\mathbb A}

\newcommand{\id}{id}

\title{On smooth curves endowed with a large automorphism $p$-group in characteristic $p>0$.}
\author{Michel Matignon and Magali Rocher.}
\date{}

\begin{document}

\maketitle

\begin{abstract}
Let $k$ be an algebraically closed field of characteristic $p>0$ and $C$ a connected nonsingular projective curve over $k$ with genus $g \geq 2$.
This paper continues the work begun in \cite{LM}, namely the study of "big actions", i.e. the pairs $(C,G)$ where $G$ is a $p$-subgroup of the $k$-automorphism group of $C$ such that$\frac{|G|}{g} >\frac{2\,p}{p-1}$. If $G_2$ denotes the second ramification group  of $G$ at the unique ramification point of the cover $C \rightarrow C/G$, we display necessary conditions on $G_2$ for $(C,G)$ to be a big action, which allows us to pursue the classification of big actions. 
Our main source of examples comes from the construction of curves with many rational points 
using ray class field theory for global function fields, as initiated by J-P. Serre and followed by \cite{Lau} and \cite{Au}. In particular, we obtain explicit examples of big actions with $G_2$ abelian of large exponent.
\end{abstract}

\section{Introduction.}
 \textit{Setting.} Let $k$ be an algebraically closed field and $C$ a connected nonsingular projective curve over $k$, with genus  $g \geq 2$.  If $char(k)=0$, Hurwitz exhibits a linear bound for the $k$-automorphism group of the curve $C$, namely: $|Aut_k(C)| \leq 84\, (g-1)$. If $char(k)=p>0$, the Hurwitz bound is no longer true when $g$ grows large, but the finiteness result still holds (cf. \cite{Sch}) and  one gets polynomial bounds on $|Aut_k(C)|$  (cf. \cite{St} and \cite{Sin}). In this situation, the full automorphism group may be very large as compared with the case $char(k)=0$. This is due to the appearance of wild ramification, which leads us to concentrate on large automorphism p-groups in $char(k)=p>0$. In this spirit,  Nakajima (cf. \cite{Na}) studies the size of Sylow $p$-subgroups of $Aut_k(C)$ and emphasizes the influence of another important invariant of the curve: the $p$-rank, denoted by $\gamma$. Indeed, if $G$ is a Sylow $p$-subgroup of $Aut_k(C)$, we deduce from \cite{Na} that $|G| \leq \frac{2\,p}{p-1} \, g$, except for $\gamma=0$. On the contrary, when $\gamma=0$, the upper bound on $|G|$ is no more linear in $g$, namely 
 $|G| \leq \max \{g, \frac{4\,p}{(p-1)^2} \, g^2\}$. As shown in \cite{St}, the
quadratic upper bound $\frac{4\,p}{(p-1)^2} \, g^2$ can really be attained, which demonstrates that, in this case, $Aut_k(C)$ may be especially large.
Following Nakajima's work, Lehr and Matignon explore  the "big actions", that is to say the pairs $(C,G)$ where $G$ is a $p$-subgroup of $Aut_k(C)$ such that
$\frac{|G|}{g}> \frac{2\,p}{p-1}$ (see \cite{LM}). In particular, they exhibit a classification of the big actions that satisfy $\frac{4}{(p-1)^2} \leq \frac{|G|}{g^2}$.

\medskip
\emph{Motivation and outline of the paper.} 
Let $(C,G)$ be a big action. As shown in \cite{LM}, there is a point of $C$, say $\infty$, such that $G$ is equal to the wild inertia subgroup $G_1$ of $G$ at $\infty$. Let $G_2$ be the second ramification group of $G$ at $\infty$ in lower notation. Then, the quotient curve $C/G_2$ is isomorphic to the projective line $\p_k^1$ and the quotient group $G/G_2$ acts as a group of translations of $\p_k^1$ fixing $\infty$, through $X \rightarrow X+y$, where $y$ runs over a subgroup $V$ of $k$. In this way, the group $G$ appears as an extension of $G_2$ by the $p$-elementary abelian group $V$ via the exact sequence:
 $$ 0 \longrightarrow G_2 \longrightarrow G=G_1  \longrightarrow V \simeq
   (\Z/\,p\, \Z)^v \longrightarrow 0$$
The aim of this paper is, on the one hand, to give necessary conditions on $G_2$ for $(C,G)$ to be a big action and, on the other hand, to display realizations of big actions with $G_2$ abelian of large exponent. In section 2, we first prove that $G_2$ must be equal to $D(G)$, the commutator subgroup of $G$. In section 3, given a big action $(C,G)$ and an additive polynomial map: $\p_k^1 \rightarrow C/G_2 \simeq \p_k^1$, we display a new big action $(\tilde{C}, \tilde{G})$ such that $\tilde{G}_2\simeq G_2$. In section 4, we demonstrate that $G_2$ cannot be cyclic except when $G_2$ has order $p$. Some of these results on $G_2$ are necessary to pursue the classification of big actions initiated by Lehr and Matignon, more precisely to explore the case: $\frac{4}{(p^2-1)^2} \leq \frac{|G|}{g^2}$. Indeed, we prove in section 5 that such an inequality requires $G_2$ to be an elementary abelian $p$-group whose order divides $p^3$. In sequel papers, M. Rocher goes further: she studies big actions with a $p$-elementary abelian $G_2$ (see \cite{MR2}), which enables her to display the classification of big actions satisfying $\frac{4}{(p^2-1)^2} \leq \frac{|G|}{g^2}$ (see \cite{MR3}). 
In section 6, following \cite{Lau} and \cite{Au}, we consider the maximal abelian extension of $K:=\F_q(X)$ ($q=p^e$) denoted by $K_S^m$, which is unramified outside $X=\infty$, completely split over the set $S$ of the finite rational places and whose conductor is smaller than $m \, \infty$, with $m \in \N$. Class field theory gives a description of the Galois group $G_S(m)$ of this extension, but also precises its upper ramification groups, which allows us to compute the genus of the extension. Moreover, it follows from the unicity and the maximality of $K_S^m$ that the group of translations $\{X \rightarrow X+y, \, y \in \F_q \}$ extends to a $p$-group of $\F_q$-automorphisms of $K_S^m$, say $G(m)$, 
with the exact sequence:
$$ 0 \longrightarrow G_S(m) \longrightarrow G(m) \longrightarrow \, \F_q \longrightarrow 0$$
This provides examples of big actions with $G_2$ abelian of exponent as large as we want, but also relates the problem of big actions to the search of algebraic curves with many rational points compared with their genera. In particular, we conclude section 6 by exhibiting specific $K$-subextensions of $K_S^m$, for a well-chosen conductor $m\,\infty$, giving examples of big actions such that $G_2 \simeq \Z/p^2\Z \times (\Z/p\Z)^t$ with a small $p$-rank, namely $t=O(log_p \,g)$. In the final section, we use Katz-Gabber theorem to highlight the link between big actions on curves and an analogous ramification condition for finite $p$-groups acting on $k((z))$.

\medskip
\textit{Notation and preliminary remarks.}
Let $k$ be an algebraically closed field of characteristic $p>0$.
We denote by $F$ the Frobenius endomorphism for a $k$-algebra. Then, $\wp$ means the Frobenius operator minus identity.
We denote by $k\{F\}$ the $k$-subspace of $k[X]$ generated by the polynomials $F^i(X)$, with $i \in \N$. It is a ring under the composition. Furthermore, for all $\alpha$ in $k$, $F\, \alpha=\alpha^p\,F$. The elements of $k\{F\}$ are the additive polynomials, i.e. the polynomials $P(X)$ of $k[X]$ such that for all $ \alpha$ and $\beta$ in $k$, $P(\alpha+ \beta) = P(\alpha)+ P(\beta)$.  Moreover, a separable polynomial is additive if and only if the set of its roots is a subgroup of $k$ (see \cite{Go} chap. 1).\\
\indent Let $f(X)$ be a polynomial of $k[X]$. Then, there is a unique polynomial $red (f)(X)$ in $k[X]$, called the reduced representative of $f$, which is $p$-power free, i.e. $red(f)(X) \in \bigoplus_{(i,p)=1} k\, X^i$, and such that $red(f)(X)=f(X)$ mod $\wp (k[X]).$ We say that the polynomial $f$ is reduced mod $\wp(k[X])$ if and only if it coincides with its reduced representative $red(f)$.
The equation $W^p-W=f(X)$ defines a $p$-cyclic \'etale cover of the affine line that we denote by $C_f$. 
Conversely, any $p$-cyclic \'etale cover of the affine line $Spec \, k[X]$ corresponds to a curve $C_f$ where $f$ is a polynomial of $k[X]$ (see \cite{Mi} III.4.12, p. 127). 
By Artin-Schreier theory, the covers $C_f$ and $C_{red(f)}$ define the same $p$-cyclic covers of the affine line.
The curve $C_f$ is irreducible if and only if $red(f) \neq 0$.\\
\indent Throughout the text, $C$ always denotes a nonsingular smooth projective curve with genus $g$  and $Aut_k(C)$ means its $k$-automorphism group.
Our main references for ramification theory are \cite{Se} and \cite{Au}.

\section{First results on "big actions".}
 To precise the background of our work, we begin by collecting and completing the first results on big actions already obtained in \cite{LM}. The expression "big actions" stands for curves endowed with a big automorphism $p$-group. The first task is to recall what we mean by "big". 

\begin{defi}
Let $G$ be a subgroup of $Aut_k (C).$
We say that the pair $(C,G)$ is a big action if $G$ is a finite $p$-group, if $g \neq 0$ and if
\begin{equation} \label{N}
\frac{ |G|}{g} > \frac{2\,p}{p-1}
\end{equation}
\end{defi}

\begin{proposition} \cite{LM}
Assume that $(C,G)$ is a big action with $g \geq 2$. Then, there is a point of $C$ (say $\infty$) such that $G$ is the wild inertia subgroup of $G$ at $\infty$: $G_1$.
Moreover, the quotient $C / G$ is isomorphic to the projective line $\p^1_k$ and the ramification locus (respectively branch locus) of the cover $\pi: \, C \rightarrow C/G$ is the 
point  $\infty$ (respectively $\pi(\infty))$.
For all $i \geq 0$, we denote by $G_i$ the $i$-th lower ramification group of $G$ at $\infty$.  
Then, 
\begin{enumerate}
\item $G_2$ is non trivial and it is strictly included in $G_1$.
\item The Hurwitz genus formula applied to $C \rightarrow C/G$ reads:
\begin{equation} \label{genus}
2 \, g = \, \sum_{i \geq 2} (|G_i|-1)
\end{equation}
\item The quotient curve $C/G_2$ is isomorphic to the projective line $\p_k^1.$ Moreover, the quotient group $G /G_2$ acts as a group of translations of the affine line $C/G_2 - \{\infty \}=Spec \,k[X]$, through $X  \rightarrow X+y$, where $y$ runs over a subgroup $V$ of $k$. Then, $V$ is an
$\F_p$-subvector space of $k$. We denote by $v$ its dimension. Thus, we obtain the exact
  sequence:
   $$ 0 \longrightarrow G_2 \longrightarrow G=G_1  \stackrel{\pi}{\longrightarrow} V \simeq
   (\Z/\,p\, \Z)^v \longrightarrow 0$$ where
 $$  \pi: 
\left\{
\begin{array}{ll}
 G \rightarrow V\\
g \rightarrow g(X)-X
\end{array}
\right.
$$ 
\item Let $H$ be a normal subgroup of $G$ such that $g_{C/H} >0$. Then, $(C/H, G/H)$ is also a big action. Moreover, the group $G/H$ fixes the image of $\infty$ in the cover $C  \rightarrow C/H$.
In particular, if $g_{C/H}=1$, then $p=2$, $C/H$ is birational to the curve $W^2+W=X^3$ and $G/H$ is isomorphic to $Q_8$, the quarternion group of order $8$ (see \cite{Si}, Appendix A, Prop. 1.2).
\end{enumerate}
\end{proposition}

\begin{remarque}
Note that, for $g=1$, one can find big actions $(C,G)$ such that $G$ is not included in a decomposition group of $Aut_k(C)$ as in Proposition 2.2.
\end{remarque}

 \noindent The following lemma generalizes and completes the last point of Proposition 2.2.
\begin{lemme}
Let $G$ a finite $p$-subgroup of $Aut_k(C)$. We assume that the quotient curve $C/G$ is isomorphic to $\p_k^1$ and that there is a point of $C$ (say $\infty$) such that $G$ is the wild inertia subgroup of $G$ at $\infty$: $G_1$.
We also assume that the ramification locus (respectively branch locus) of the cover $\pi: \, C \rightarrow C/G$ is the 
point  $\infty$ (respectively $\pi(\infty))$. Let $G_2$ be the second ramification group of $G$ at $\infty$ and $H$ a subgroup of $G$. 
\begin{enumerate}
\item Then, $C/H$ is isomorphic to $\p_k^1$ if and only if $H \supset G_2$. 
\item In particular, if $(C,G)$ is a big action with $g \geq 2$ and if
$H$ is a normal subgroup of $G$ such that $H \subsetneq G_2$, then $g_{C/H}>0$ and $(C/H,G/H)$ is also a big action. 
\end{enumerate}
\end{lemme}

\noindent \textbf{Proof}:
When applied to the cover $C \rightarrow C/G \simeq \p_k^1$, the Hurwitz genus formula (see e.g. \cite{St93}) reads:
$2(g-1)=2|G|\, (g_{C/G}-1) + \sum_{i\geq 0}\, (|G_i|-1)$.
When applied to the cover $C \rightarrow C/H$, it yields:
$2(g-1)=2|H|\, (g_{C/H}-1) + \sum_{i\geq 0} \,(|H \cap G_i|-1)$.
Since $H \subset G=G_0=G_1$, it follows that:
$$ 2 |H| g_{C/H}= -\, 2(|G|-|H|) + \sum_{i\geq 0} \, (|G_i|-|H \cap G_i|)= \sum_{i\geq 2} \,\, (|G_i|-|H \cap G_i|) $$
Therefore, $g_{C/H}=0$ if and only if for all $i \geq 2$, $G_i=H \cap G_i$, i.e. $G_i \subset H$, which is equivalent to $G_2 \subset H$. The second point then derives from Proposition 2.2.4.
$\square$
\medskip

The very first step to study big actions is to precise their description when $G_2 \simeq \Z/p \Z$.
The following proposition aims at gathering and reformulating the results already obtained for this case in \cite{LM} (cf. Prop. 5.5, 8.1 and 8.3).

\begin{proposition} \cite{LM}.
Let $(C,G)$ be a big action, with $g \geq 2$, such that $G_2 \simeq \Z/p \Z$. 
\begin{enumerate}
\item Then, $C$ is birational to the curve $C_f: \, W^p-W=f(X)=X\,S(X)+c\, X \in k[X]$, where $S$ in $k\{F \}$ is an additive polynomial with degree $s\geq 1$ in $F$. If we denote by $m$ the degree of $f$, then $m=1+p^s=i_0$, where $i_0 \geq 2$ is the integer such that:
$$G=G_0=G_1 \supsetneq G_2 = G_3=\cdots=G_{i_0} \supsetneq G_{i_0+1} =\cdots=$$
\item Write $ S(F) = \sum _{j=0}^{s} a_j F^j$, with $a_{s} \neq 0$.
Then, following \cite{El} (section 4), we can define an additive polynomial related to $f$, called the "palindromic polynomial" of $f$: $$Ad_{f}:= \frac{1}{a_{s}} \, F^{s} \,  (\sum _{j=0}^ {s} \,a_j\, F^j + F^{-j}\, a_j)$$
 The set of roots of $Ad_f$, denoted by $Z(Ad_f)$, is an $\F_p$-subvector space of $k$,
isomorphic to $(\Z/p\Z)^{2s}$. 
 Besides, $Z(Ad_f)= \{y \in k,\, f(X+y)-f(X)=0 \, \mod \, \wp(k[X] \}$.
\item Let $G_{\infty,1}$ be the wild inertia subgroup of $Aut_k(C)$ at $\infty$.
Then, $G_{\infty,1}$  is a central extension of $ \Z/p \Z$ by the elementary abelian $p$-group $Z(Ad_f)$ which can be identified with a subgroup of translations $ \{ X \rightarrow X+y ,\, \,  y \in k \}$  of the affine line. 
Furthermore, if we denote by $Z(G_{\infty,1})$ the center of $G_{\infty,1}$ and
by $D(G_{\infty,1})$ its commutator subgroup, $Z(G_{\infty,1})=D(G_{\infty,1})=<\sigma>$, where $\sigma(X)=X$ and $\sigma(W)=W+1$. 
Thus, we get the following exact sequence:
$$ 0 \longrightarrow Z(G_{\infty,1})=D(G_{\infty,1})\simeq \Z/p \Z \longrightarrow G_{\infty,1}  \stackrel{\pi}{\longrightarrow}  Z(Ad_f) \simeq (\Z / p \Z) ^{2s} \longrightarrow 0$$ where  
 $$  \pi: 
\left\{
\begin{array}{lc}
 G_{\infty,1} \rightarrow Z(Ad_f) \simeq (\Z / p \Z) ^{2s}\\
g \rightarrow g(X)-X
\end{array}
\right.
$$ 
For $p>2$,  $G_{\infty,1}$ is the unique extraspecial group with exponent $p$ and order $p^{2s+1}$. The case $p=2$ is more complicated (see \cite{LM} 4.1).

\item There exists an $\F_p$-vector space $V \subset Z(Ad_f) \simeq (\Z / p \Z) ^{2s}$ such that
$G=\pi^{-1}(V) \subset G_{\infty,1}$ and such that we get the exact sequence:
$$ 0 \longrightarrow G_2 \simeq \Z/p \Z \longrightarrow G \stackrel{\pi}{\longrightarrow} V \longrightarrow 0.$$
\end{enumerate}
\end{proposition}

Therefore, the key idea to study big actions is to use Proposition 2.2.4 and Lemma 2.4.2 to go back to the well-known situation described above. This motivates the following 

\begin{theoreme}
Let $(C,G)$ be a big action with $g \geq 2$.
Let $\mathcal{G}$ be a normal subgroup in $G$ such that $\mathcal{G}$ is strictly included in $G_2$. Then, there exists a group $H$, normal in $G$, such that $\mathcal{G} \subset H \subsetneq G_2$ and $[G_2:H]=p$. In this case, $(C/H,G/H)$ enjoys the following properties. 
 \begin{enumerate}
\item The pair $(C/H,G/H)$ is a big action and the exact sequence of Proposition 2.2:
 $$ 0 \longrightarrow G_2  \longrightarrow G \stackrel{\pi}{\longrightarrow}  V \longrightarrow 0$$ induces the following one:
$$ 0 \longrightarrow G_2/H=(G/H)_2 \simeq \Z/p\Z \longrightarrow G/H \stackrel{\pi}{\longrightarrow}  V \longrightarrow 0$$
\item The curve $C/H$ is birational to $C_f$: $W^p-W=f(X)=X \,S(X)+c\,X \in k[X]$, where $S$ is an additive polynomial of degree $s \geq 1$ in $F$. Let $Ad_f$ be the palindromic polynomial related to $f$ as defined in Proposition 2.5. Then, $V \subset Z(Ad_f)\simeq (\Z/p\Z)^{2s}$.
\item Let $E$ be the wild inertia subgroup of $Aut_k(C/H)$ at $\infty$. We denote by $D(E)$ its commutator subgroup of $E$ and by $Z(E)$ its center.
Then, $E$ is an extraspecial group of order $p^{2s+1}$ and
 $$ 0 \longrightarrow D(E)=Z(E) \simeq \Z/p \Z \longrightarrow E \stackrel{\pi}{\longrightarrow}  Z(Ad_f) \simeq (\Z
 /p \Z) ^{2s} \longrightarrow 0$$
 \item Moreover, $G/H$ is a normal subgroup in $E$.
 It follows that $G_2$ is equal to $D(G)$, the commutator subgroup of $G$, which is also equal to $D(G)G^p$, the Frattini subgroup of $G$.
 \end{enumerate}
\end{theoreme}

\noindent \textbf{Proof:} First of all, the existence of the group $H$ comes from \cite{Su1} (Chap. 2, Thm. 1.12). We deduce from Lemma 2.4.2 that $(C/H, G/H)$ is still a big action. Then, $G=G_1 \supsetneq G_2$ (resp. $G/H=(G/H)_1 \supsetneq (G/H)_2$). As the first jump always coincides in
lower and upper ramification, it follows that $G_2=G^2$ (resp.$(G/H)_2=(G/H)^2$). By \cite{Se} (Second Part, Chap. IV, Prop. 14), $(G/H)_2=(G/H)^2= G^2H/H=G_2H/H=G_2/H$. The first assertion follows.
The second and the third point directly derive from Proposition 2.5. \\
\indent We now prove the last statement. By Proposition 2.5, $Z(E)=(G/H)_2=G_2/H \subset G/H$. So, $G/H$ is a subgroup of $E$ containing $Z(E)$. Moreover, since $(\Z/p\Z)^{2s}$ is abelian, $\pi(G/H)$ is normal in $E/Z(E)$. It follows that $G/H$ is normal in $E$. We eventually show that $G_2=D(G)$.  On the one hand, since $G/G_2$ is abelian, $D(G)$ is included in $G_2.$ On the other hand, assume that $D(G)$ is strictly included in $G_2$. Then, the first point applied to $\mathcal{G}=D(G)$ ensures the existence of a group $H$, normal in $G$, with $D(G) \subset H \subset G_2$, $[G_2:H]=p$ and such that $(C/H, G/H)$ is a big action. Since $D(G) \subset H$, $G/H$ is an abelian subgroup of $E$. As  $G/H$ is also a normal group in $E$, \cite{Hu} (Satz 13.7) implies
 $|G/H| \leq p^{s+1}$. Hence $\frac{|G/H|}{g_{C/H}} \leq \frac{2\,p}{p-1}$, which contradicts condition \eqref{N} for the big action $(C/H, G/H)$. It follows that $D(G)=G_2$. In addition, as $G/G_2$ is an elementary abelian $p$-group, then $G^p=G_1^p \subset G_2=D(G)$. As a consequence, $G_2=D(G)G^p$ which is equal to the Frattini subgroup of $G$, since $G$ is a $p$-group. $\square$

\begin{remarque}
When applying Theorem 2.6 to $\mathcal{G}=G_{i_0+1}$, where $i_0$ is defined as in 
Proposition 2.5, one obtains Theorem 8.6(i) of \cite{LM}. In particular, for all big actions $(C,G)$ with $g \geq 2$, there exists an index $p$-subgroup $H$ of $G_2$, normal in $G$, such that $(C/H,G/H)$ is a big action with $C/H$ birational to $W^p-W=f(X)=X\,S(X)+c\,X \in k[X]$, where $S$ is an additive polynomial of degree $s\geq 1$ in $F$. Note that, in this case, $i_0=1+p^s$.
\end{remarque}

As $G_2$ cannot be trivial for a big action, we gather from the last point of Theorem 2.6 the following result.

\begin{corollaire}
Let $(C,G)$ be a big action with $g \geq 2$. Then $G$ cannot be  abelian.
\end{corollaire}

It is natural to wonder whether $G_2$ can be non abelian. Although we do not know yet the answer to this question, we can mention a special case in which $G_2$ is always abelian, namely:

\begin{corollaire}
Let $(C,G)$ be a big action with $g \geq 2$. If the order of $G_2$ divides $ p^3$, then $G_2$ is abelian.
\end{corollaire}

\noindent \textbf{Proof:} There is actually only one case to study, namely: $|G_2|=p^3$.
We denote by $Z(G_2)$ the center of $G_2$.
The case $|Z(G_2)| = 1$ is impossible since $G_2$ is a $p$-group. 
If $|Z(G_2)| = p$, then $Z(G_2)$ is cyclic. But, as $G_2$ is a $p$-group, normal in $G$ and included in $D(G)$ (see Theorem 2.6), \cite{Su2} (Prop. 4.21, p. 75) implies that $G_2$ is also cyclic, which contradicts the strict inclusion of $Z(G_2)$ in $G_2$. 
If $|Z(G_2)| = p^2$, then $G_2/Z(G_2)$ is cyclic and $G_2$ is abelian, which leads to the same contradiction as above. This leaves only one possibility: $|Z(G_2)|= p^3$, which means that $G_2=Z(G_2)$. $\square$

\begin{corollaire}
Let $(C,G)$ be a big action with $g \geq 2$. We keep the notation of Proposition 2.2. Let  $G_{\infty,1} $ be the wild inertia subgroup of $Aut_k(C)$ at $\infty$.
Then, $(C,G_{\infty,1} )$ is a big action whose second lower ramification group is equal to $D(G_{\infty,1})=D(G)$. In particular, $G$ is equal to $G_{\infty,1}$
if and only if $|G/D(G)| =|G_{\infty,1} /D(G_{\infty,1})|$.
\end{corollaire}

\noindent \textbf{Proof:} As $G$ is included in $G_{\infty,1}$, then $D(G) \subset D(G_{\infty,1})$.
If the inclusion is strict, one can find a subgroup $\mathcal{G}$ such that $G \subsetneq \mathcal{G} \subset G_{\infty,1}$ with $[\mathcal{G}:G]=p$ (see \cite{Su1}, Chap. 2, Thm. 19). Then, $G$ is a normal subgroup of $\mathcal{G}$.
It follows that $D(G)$ is also a normal subgroup of $\mathcal{G}$. 
As $|G| \leq |\mathcal{G}|$, the pair $(C,\mathcal{G})$ is a big action. So, by Theorem 2.6, $\mathcal{G}_2=D(\mathcal{G})$. Since $D(G)$ is normal in $\mathcal{G}$ and $g(C/D(G))=0$, we gather from Lemma 2.4.1 that $D(G)=\mathcal{G}_2=D(\mathcal{G})$. The claim follows by reiterating the process. $\square$

\begin{remarque}
Let $(C,G_{\infty,1})$ be a big action as in Corollary 2.10. Then, $G_{\infty,1}$ is a $p$-Sylow subgroup of $Aut_k(C)$. Moreover, we deduce from \cite{Gi} (Thm. 1.3) that $G_{\infty,1}$ is the unique $p$-Sylow subgroup of $Aut_k(C)$ except in four special cases: the hyperelliptic curves: $W^{p^n}-W=X^2$ with $p>2$, the Hermitian curves and the Deligne-Lusztig curves arising from the Suzuki groups and the Ree groups (see equations in \cite{Gi}, Thm. 1.1).
\end{remarque}

\section{Base change and big actions.}

Starting from a given big action $(C,G)$, we now display a way to produce a new one: $(\tilde{C}, \tilde{G})$, with $\tilde{G}_2\simeq G_2$ and $g_{\tilde{C}}=p^{s} \, g_C$. The main tool is a base change associated with an additive polynomial map: $\p_k^1  \stackrel{S}{\longrightarrow} C/G_2\simeq \p_k^1 $.

\begin{proposition}
Let $(C,G)$ be a big action with $g \geq 2$.
We denote by $L:=k(C)$ the function field of the curve $C$, by $k(X):=L^{G_2}$ the subfield of $L$ fixed by $G_2$  and by $k(T):=L^{G_1}$, with $T=\prod_{v \in V} (X-v)$.
Write $X=S(Z)$, where $S(Z)$ is a separable additive polynomial of $k[Z]$ with degree $p^{s}$, $s \in \N$.
\begin{enumerate}
\item Then, $L$ and $k(Z)$ are linearly disjoined over $k(X)$.
\item Let $\tilde{C}$ be the smooth projective curve over $k$ with function field $k(\tilde{C}):=L[Z]$. Then, $k(\tilde{C})/k(T)$ is a Galois extension with group $\tilde{G} \simeq G \times (\Z/p\Z)^{s}$. Furthermore, $g_{\tilde{C}}=p^{s} \, g_C$. It follows that $\frac{|\tilde{G}|}{g_{\tilde{C}}}= \frac{|G|}{g}$. So, $(\tilde{C},\tilde{G})$ is still a big action with second ramification group $\tilde{G_2}\simeq G_2\times \{0\} \subset G \times (\Z/p\Z)^{s}$. This can be illustrated by the following diagram:
\end{enumerate}
$$
\begin{array}{clc}
C &   \longleftarrow & \tilde{C}\\
\downarrow  & & \downarrow  \\
  C/G_2\simeq \p_k^1&  \stackrel{S}{\longleftarrow}& \p_k^1\\
\end{array}
$$
\end{proposition}

\noindent The proof of this proposition requires two preliminary lemmas.

\begin{lemme}
Let $K:=k((z))$ be a formal power series field over $k$. 
Let $K_1/K$ be a Galois extension whose group $\mathcal{G}$ is a $p$-group. Let $K_0/K$ be a $p$-cyclic extension. Assume that $K_0$ and $K_1$ are linearly disjoined over $K$. Put $L:=K_0K_1$.
$$
\begin{array}{clc}
K_1 & - & L=K_0 K_1\\
\quad \\
 \mathcal{G} \, \vert& & \vert  \quad \\
\quad \\
K & -& K_0
\end{array}
$$
We suppose that the conductor of $K_0/K$ (see e.g. \cite{Se} Chap. 15, Cor. 2) is $2$. Then, $L/K_1$ also has conductor $2$.
\end{lemme}

\noindent \textbf{Proof:} 
Consider a chief series of $\mathcal{G}$ (cf. \cite{Su1}, Chap. 2, Thm. 1.12), i.e. a series such that: $$\mathcal{G} =\mathcal{G}_0 \supsetneq \mathcal{G}_1 \cdots \supsetneq \mathcal{G}_n=\{0\}$$
with $\mathcal{G}_i$ normal in $\mathcal{G}$ and $[\mathcal{G}_{i-1}: \mathcal{G}_i]=p$.
One shows, by induction on $i$, that the conductor of the extension $K_0K_1^{ \mathcal{G}_i}/ K_1^{ \mathcal{G}_i}$ is $2$. Therefore, one can assume $\mathcal{G} \simeq \Z/p\Z$.
In this case, $L/k((z))$ is a Galois extension with group 
$G \simeq (\Z/p\Z)^2$. Write the ramification filtration of $G$ in lower notation:
$$G=G_0=\cdots=G_{i_0} \supsetneq G_{i_0+1} =\cdots$$
\begin{enumerate}
\item First, assume that $G_{i_0+1}=\{0\}$. 
Then, an exercise shows that, for any index $p$-subgroup $H$ of $G$, the extensions $L/L^H$ (case $(\alpha)$) and $L^H/K$ (case $(\beta)$) are $p$-cyclic with conductor $i_0+1$.
When applied to $H=Gal(L/K_0)$, case $(\beta)$ gives $i_0=1$. Therefore, one concludes by applying case $(\alpha)$ to $H=Gal(L/K_1)$.
\item Now, assume that $G_{i_0+1}\neq \{0\}$. As above, let $H$ be an index $p$-subgroup of $G$. An exercise using the classical properties of ramification theory shows that:
\begin{enumerate}
\item If $H=G_{i_0+1}$, then $L/L^{H}$ (resp. $L^{H}/K$) is a $p$-cyclic extension with conductor $i_0+i_1+1$ (resp. $i_0+1$).
\item If $H\neq G_{i_0+1}$, then $L/L^{H}$ (resp. $L^{H}/K$) is a $p$-cyclic extension with conductor $i_0+1$ (resp. $i_0+\frac{i_1}{p}+1$).
\end{enumerate} 
Apply this result to $H:=Gal(L/K_0)$. Since $K_0/K$ has conductor $2$, it follows that $i_0+1=2$, so $i_0=1$  and $Gal(L/K_0)=G_{i_0+1}$. Therefore, $Gal(L/K_1) \neq G_{i_0+1}$ and we infer from case (b) that $L/K_1$ has conductor $i_0+1=2$. $\square$
\end{enumerate}

\begin{lemme}
Let $W$ be a finite $\F_p$-subvector space of $k$. 
Let $W_1$ and $W_2$ be two $\F_p$-subvectors spaces of $W$ such that $W =W_1 \bigoplus W_2$.
Define $T:=\prod_{w \in W} (Z-w)$ and $T_i:=\prod_{w \in W_i} (Z-w)$, for $i$ in $\{1,2\}$. Then, $k(T) \subset k(T_i) \subset k(Z)$. Moreover,
\begin{enumerate}
\item  The extensions $k(T_1)/k(T)$ and $k(T_2)/k(T)$ are linearly disjoined over $k(T)$.
\item For all $i$ in $\{1,2\}$, $k(Z)/k(T)$ (resp. $k(Z)/k(T_i)$) is a Galois extension with group isomorphic to $W$ (resp. $W_i$).
 \item For all $i$ in $\{1,2\}$, $k(T_i)/k(T)$ is a Galois extension with group isomorphic to $\frac{W}{W_i}$.
\end{enumerate}
This induces the following diagram:
$$
\begin{array}{lll}
k(T_1) &    \stackrel{W_1}{-} & k(Z)\\
& &\quad \\
 \vert \, \frac{W}{W_1} & & \vert \, W_2  \\
& &\quad \\
k(T) &   \stackrel{\frac{W}{W_2}}{-} & k(T_2)
\end{array}
$$
\end{lemme}

\noindent \textbf{Proof:}
Use for example \cite{Go} (1.8). $\square$

\medskip

\noindent \textbf{Proof of Proposition 3.1:}
\begin{enumerate}
\item The first point derives from Lemma 2.4.1.
\item Put $W:=S^{-1}(V)$, with $V$ defined as in Proposition 2.2.3, $W_1:=S^{-1}(\{0\})\simeq (\Z/p\Z)^{s}$, since $S$ is an additive separable polynomial of $k[Z]$ with degree $p^{s}$ (see e.g. \cite{Go} chap. 1). Call $W_2$ any $\F_p$-subvector space of $W$ such that $W=W_1 \bigoplus W_2$. Then, Lemma 3.3 applied to the extension $k(Z)/k(T)$ induces the following diagram:
\medskip
$$
\begin{array}{clc}
L=k(C) &   - & k(\tilde{C})\\
\quad \\
G_2 \,\vert  & & \vert  \\
\quad \\
L^{G_2}=k(X)=k(Z)^{W_1} &  \stackrel{W_1}{-} & k(Z)\\
\quad \\
 \frac{W}{W_1} \, \vert & & \vert \, W_2  \\
\quad \\
L^{G_1}=k(T)=k(Z)^{W} & \stackrel{\frac{W}{W_2}}{-}  & k(Z)^{W_2}
\end{array}
$$
\medskip
In particular, Lemma 3.3  implies that $k(Z)^{W_1} \cap k(Z)^{W_2}=k(T)$. 
Since $k(C)\cap k(Z)=k(X)$ (cf. first point of the proposition), we deduce that $k(C)$ and $k(Z)^{W_2}$ are linearly disjoined over $k(T)$. As $k(Z)^{W_2}/k(T)$ is a Galois extension with group $\frac{W}{W_2} \simeq W_1 \simeq (\Z/p \Z)^{s}$, it follows that $k(\tilde{C})/k(T)$ is a Galois extension with group $\tilde{G} \simeq Gal(k(C)/k(T)) \times Gal(k(Z)^{W_2} /k(T)) \simeq G \times (\Z/p\Z)^{s}$.\\
\indent Now, consider a flag of $\F_p$-subvector spaces of $W_1$:
$$W_1=W_1^{(1)} \supsetneq  W_1^{(2)} \supsetneq  \cdots \supsetneq  W_1^{(s+1)}=\{0\}$$
such that $[W_1^{(i-1)}:W_1^{(i)}]=p$. It induces the following inclusions:
$$k(Z)=k(Z)^{W_1^{(s+1)}} \supsetneq k(Z)^{W_1^{(s)}} \supsetneq \cdots \supsetneq k(Z)^{W_1^{(1)}} =k(X)$$
Then, apply Lemma 3.2 to $K_1/K$: the completion at $\infty$ of the extension $k(C)/k(X)$, whose group $G_2$ is a $p$-group, and to $K_0/K$: the completion at $\infty$ of the $p$-cyclic extension $k(Z)^{W_1^{(i)}}/k(Z)^{W_1^{(i-1)}}$ whose conductor is $2$. By induction, we thus prove that the extension $k(\tilde{C})/k(T)$ also has conductor $2$. It follows from the Hurwitz genus formula that $g_{\tilde{C}}=p^{s} \, g_C$. Finally, the last statement on $\tilde{G}_2$ is a consequence of Lemma 2.4.1. $\square$
\end{enumerate}

\begin{remarque}
Under the conditions of Proposition 3.1, it can happen that $G$ is a $p$-Sylow subgroup of $Aut_k(C)$ without
$\tilde{G}$ being a $p$-Sylow subgroup of $Aut_k(\tilde{C})$. \\
\indent Indeed, take $C$ : $W^p-W=X^{1+p}$ and $S(Z)=Z^p-Z$. Then, $\tilde{C}$ is parametrized by $\tilde{W}^p-\tilde{W}=(Z^p-Z)\, (Z^{p^2}-Z^p)=-Z^2+2\, Z^{1+p} -Z^{1+p^2}$ mod $\wp(k[Z])$. We denote by $G_{\infty,1}(C)$ (resp. $G_{\infty,1}(\tilde{C})$) the 
wild inertia subgroup of $Aut_k(C)$ (resp. $Aut_k(\tilde{C})$) at $X=\infty$ (resp. $Z=\infty$). Note that $G_{\infty,1}(C)$ (resp. $G_{\infty,1}(\tilde{C})$) is a $p$-Sylow subgroup of $Aut_k(C)$ (resp. $Aut_k(\tilde{C})$). Take $G:=G_{\infty,1}(C)$.
From Proposition 2.5, we deduce that $|\tilde{G}|=p\,|G|=p\,|G_{\infty,1}(C)|=p^4$, whereas $|G_{\infty,1}(\tilde{C})|=p^5$.
\end{remarque}

\section{A new step towards a classification of big actions.}

 If big actions are defined through the value taken by the quotient $\frac{|G|}{g}$, it occurs that the key criterion to classify them is the value of another quotient: $\frac{|G|}{g^2}$. Indeed, the quotient $\frac{|G|}{g^2}$ has, to some extent, a "sieve" effect among big actions. In what follows, we pursue the work of Lehr and Matignon who describe big actions for the two highest possible values of this quotient, namely $\frac{|G|}{g^2}=\frac{4\, p}{(p-1)^2}$ and $\frac{|G|}{g^2}=\frac{4}{(p-1)^2}$ (cf. \cite{LM} Thm. 8.6). More precisely, we investigate the big actions $(C,G)$ that satisfy:
\begin{equation} \label{eq*}
M:= \frac{4}{(p^2-1)^2} \leq \frac{|G|}{g^2} 
\end{equation}
The choice of the lower bound $M$ can be explained as follows: as shown in the proof of (\cite{LM}, Thm. 8.6), a lower bound $M$ on the quotient $\frac{|G|}{g^2}$ involves an upper bound on the order of the second ramification group, namely:
\begin{equation} \label{eq**}
|G_2| \leq  \frac{4}{M}  \frac{|G_2/ G_{i_0+1}| ^2}{( |G_2/ G_{i_0+1}| -1)^2} 
\end{equation}
where $i_0$ is defined as in Proposition 2.5. Therefore, we have to choose $M$ small enough to obtain a wide range of possibilities for the quotient, but meanwhile large enough to get serious restrictions on the order of $G_2$. The optimal bound seems to be $M:=\frac{4}{(p^2-1)^2}$, insofar as, for such a choice of M, the upper bound on $G_2$ implies that its order divides $p^3$, and then that $G_2$ is abelian (cf. Corollary 2.9).

\begin{proposition}
Let $(C,G)$ be a big action with $g \geq 2$ satisfying condition \eqref{eq*}. Then, the order of $G_2$ divides $p^3$. It follows that $G_2$ is abelian.
\end{proposition}

\noindent \textbf{Proof:}
Put $p^m:=|G_2/G_{i_0+1}|$, with $m \geq 1$, and $A_m := \frac{4}{M}  \frac{|G_2/ G_{i_0+1}|}{(|G_2/ G_{i_0+1}| -1)^2}=\frac{4}{M}
 \frac{p^{m}}{(p^{m}-1)^2}$. Then, inequality \eqref{eq**} reads: $1< |G_2| =p^{m} |G_{i_0+1}|\leq  p^{m} A_m$, which gives: $1 \leq |G_{i_0+1}| \leq A_m$.
Since $(A_m)_{m \geq 1}$ is a decreasing sequence with $A_4 <1$, we conclude that $m \in \{1,2,3 \}$.\\ 
\indent If $m=3$, then $1 \leq |G_{i_0+1}| \leq A_3 <p$. So $|G_{i_0+1}|=1$ and $|G_2|=p^3$. 
If $m=2$, then $1 \leq |G_{i_0+1}| \leq A_{2}=p^2$.  So $|G_2| = p^2 \, |G_{i_0+1}|$, with $ |G_{i_0+1}| \in \{1,p,p^2 \}$. This leaves only one case to exclude, namely $|G_{i_0+1}|=p^2$.
In this case, $|G_2|=p^4$ and formula \eqref{genus} yields a lower bound on the genus, namely:
$ 2\, g \geq  \, (i_0-1) (p^4-1).$
 Let $s$ be the integer defined in Remark 2.7. Then, $i_0=1+p^s$. Besides, by Theorem 2.6, $V \subset (\Z/p\Z)^{2s}$. Consequently, $|G|=|G_2||V| \leq p^{4+2s}$ and  $$\frac{|G|}{g^2} \leq \frac{4\, p^{4+2s}}{p^{2s}(p^4-1)^2}=\frac{4}{(p^2-1)^2} \frac{p^4}{(p^2+1)^2} <\frac{4}{(p^2-1)^2}$$
which contradicts equality \eqref{eq*}. \\
\indent If $m=1$, then $1 \leq |G_{i_0+1}| \leq A_{1}$ with
 $ A_{1}:=p \, (p+1)^2 <
\left\{
\begin{array}{ll}
p^4 \,, \quad if \,  p \geq 3 \\
p^5 \,, \quad if \, p=2
\end{array}
\right.
$.\\ 
Since $G_{i_0+1}$ is a $p$-group, we get:
$
\left\{
\begin{array}{ll}
1 \leq |G_{i_0+1}| \leq p^3 \,, \quad if \,  p \geq 3 \\
1 \leq |G_{i_0+1}| \leq p^4 \,, \quad if \, p=2
\end{array}
\right.
$.
As $|G_2|=p\,|G_{i_0+1}|$, there are two cases to exclude: $|G_{i_0+1}|=p^{3+\epsilon}$, with $\epsilon =0$ if $p \geq 3$ and $\epsilon \in \{0,1\}$ if $p=2$. Then $|G_2|=p^{4+\epsilon}$. If $\epsilon =0$, we are in the same situation as in the previous case. If $\epsilon=1$, \eqref{genus} yields $2\, g \geq (i_0-1) (p^{5}-1).$
Since this case only occurs for $p=2$, we eventually get an inequality: 
$$\frac{|G|}{g^2} \leq \frac{4\, p^{5+2s}}{p^{2s}\,(p^5-1)^2}=  \frac{128}{961} <\frac{4}{9}= \frac{4}{(p^2-1)^2}$$
which contradicts condition \eqref{eq*}. Therefore, the order of $G_2$ divides $p^3$.
Then, we gather from Corollary 2.9 that $G_2$ is abelian. $\square$
\medskip

But we can even prove better: under these conditions, $G_2$ has exponent $p$. 
\begin{proposition} 
Let $(C,G)$ be a big action with $g \geq 2$ satisfying condition \eqref{eq*}.Then $G_2$ is abelian with exponent $p$.
\end{proposition}

\noindent \textbf{Proof:} By Proposition 4.1, $G_2$ is abelian, with order dividing $p^3$.
As a consequence, if $G_2$ has exponent strictly greater than $p$, either $G_2$ is cyclic with order $p^2$ or $p^3$, or $G_2$ is isomorphic to $\Z/p^2 \Z \times \Z/p \Z.$
We begin with a lemma excluding the second case. Note that one can find big actions $(C,G)$ with $G_2$ abelian of exponent $p^2$. Nevertheless, it requires the $p$-rank of $G_2$
 to be large enough (see section 6).

\begin{lemme}
Let $(C,G)$ be a big action with $g \geq 2$ satisfying condition \eqref{eq*}. Then $G_2$ cannot be isomorphic to $\Z/p^2 \Z \times \Z/p \Z$.
\end{lemme}

\noindent \textbf{Proof:}
Assume $G_2 \simeq \Z/p^2 \Z \times \Z/p \Z$. Then, the lower ramification filtration of $G$ reads as in one of the four following cases:

\begin{description}
\item i) $G=G_1 \supsetneq G_2 \simeq \Z/p^2 \Z \times \Z/p\Z \supset  G_{i_0+1} \simeq \Z/p\Z \supset G_{i_0+i_1+1} = \{0 \}.$
\item ii)  $G=G_1 \supsetneq G_2 \simeq \Z/p^2 \Z \times \Z/p\Z \supset  G_{i_0+1}\simeq (\Z/p\Z)^2 \supset G_{i_0+i_1+1} = \{0 \}.$
\item iii)  $G=G_1 \supsetneq G_2 \simeq \Z/p^2 \Z \times \Z/p\Z \supset  G_{i_0+1} \simeq (\Z/p\Z)^2 \supset G_{i_0+i_1+1} \simeq \Z/p\Z \supset  G_{i_0+i_1+i_2} = \{0 \}.$
\item iv) $G=G_1 \supsetneq G_2 \simeq \Z/p^2 \Z \times \Z /p \Z \supset  G_{i_0+1} \simeq \Z/p^2 \Z \supset G_{i_0+i_1+1} \simeq \Z/p\Z \supset  G_{i_0+i_1+i_2} = \{0 \}.$
\end{description}

We now focus on the ramification filtration of $G_2$, temporary denoted by $H$ for convenience. Then, for all $i \geq 0$, the lower ramification groups of $H$ are: $H_i=H \cap G_i$. \\
In case i), the lower ramification of $H$ reads:
 $$H=H_0=\cdots=H_{i_0} \simeq \Z/p^2 \Z \times \Z/p\Z \supset  H_{i_0+1}=\cdots=H_{i_0+i_1} \simeq \Z/p\Z \supset H_{i_0+i_1+1} = \{0 \}.$$
Consider the upper ramification groups:
$H^{\nu_0}=H^{\varphi(i_0)}=H_{i_0}$ and $H^{\nu_1}=H^{\varphi(i_0+i_1)}=H_{i_0+i_1}$, 
 where $\varphi$ denotes the Herbrand function. Then, the ramification filtration in upper notation reads:
$$H^0=\cdots=H^{\nu_0} \simeq \Z /p^2 \Z \times \Z/p\Z \supset  H^{\nu_0+1}=\cdots
=H^{\nu_1} \simeq \Z/p\Z \supset H^{\nu_1+1} =\{0 \}.$$  
Since $H$ is abelian, it follows from Hasse-Arf theorem that $\nu_0$ and $\nu_1$ are integers. Consequently, the formula:
$$
\forall \, m \in \N, \quad  
\varphi(m)+1= \frac{1}{|H_0|} \, \sum_{i=0}^m |H_i|
$$ gives $\nu_0=i_0$ and $\nu_1 = i_0+\frac{i_1}{p^2}$. 
 Besides, \cite{Ma} (Thm. 6) implies $H^{\nu_0} \supsetneq H^{p\,\nu_0} \supset (H^{\nu_0})^p$ with  $ (H^{\nu_0})^p=H^p=G_2^p \simeq \Z/p \Z$. Thus, $ H^{p\nu_0}  \supset H^{\nu_1}$, which involves: $ p \nu_0  \leq \nu_1$ and $i_1  \geq p^2(p-1)i_0$.
Then, the Hurwitz genus formula applied to $C \rightarrow C/H \simeq \p^1_k$ yields a lower bound for the genus:
$$2\, g =(i_0-1) (|H|-1)+i_1(|H_{i_0+1}|-1) \geq (p-1) (i_0+1) (p^3+p+1).$$
Let $s$ be the integer defined in Remark 2.7. Then, $i_0=1+p^s$.
Moreover, by Theorem 2.6, $|G|=|G_2||V| \leq  p^{3+2s}$. It follows that $ \frac{|G|}{g^2} \leq \frac{4}{(p^2-1)^2} \frac{p^3(p+1)^2}{(p^3+p+1)^2} $. 
Since $\frac{p^3(p+1)^2}{(p^3+p+1)^2}<1$ for $p \geq 2$, this contradicts condition \eqref{eq*}. 
\medskip

In case ii), the lower ramification filtration of $H$ reads:
$$H=H_0=\cdots=H_{i_0} \simeq \Z/p^2 \Z \times \Z/p\Z \supset  H_{i_0+1}=\cdots H_{i_0+i_1} \simeq (\Z/p\Z)^2 \supset H_{i_0+i_1+1} =\{0 \}.$$ 
 Keeping the same notation as in case i), the upper ramification filtration reads:
$$H=H^0=\cdots=H^{\nu_0} \simeq \Z /p^2 \Z \times \Z/p\Z \supset  H^{\nu_0+1}=\cdots
=H^{\nu_1}\simeq (\Z/p\Z)^2 \supset H^{\nu_1+1} =\{0 \}.$$ 
with
$\nu_0=\varphi(i_0)=i_0$ and $\nu_1=\varphi(i_0+i_1)=i_0+\frac{i_1}{p}$.
Once again, $H^{p\nu_0} \supset (H^{\nu_0})^p  \simeq \Z/p \Z$ implies $ H^{p \,\nu_0} \supset H^{\nu_1}$, which involves $p \, \nu_0  \leq \nu_1$ and  $i_1 \geq i_0 \, p \, (p-1).$
Then, the Hurwitz genus formula yields: $$2\, g=(i_0-1) (|H|-1)+i_1(|H_{i_0+1}|-1) \geq (p-1) \, p^s \, (p^3+p^2+1) \geq (p-1)p^s (p^3+p+1).$$ Thus, we get the same lower bound on the genus as in the preceding case, hence the same contradiction. 
\medskip

In case iii), the lower ramification filtration of $H$ reads:
$$H_{i_0} \simeq \Z/p^2 \Z \times \Z/p\Z \supset  H_{i_0+1}=\cdots=H_{i_0+i_1} \simeq (\Z/p\Z)^2 \supset H_{i_0+i_1+1} =\cdots=H_{i_0+i_1+i_2} \simeq \Z / p \Z \supset \{ 0 \}. $$ 
Keeping the same notation as above and introducing  $H^{\nu_2}=H^{\varphi(i_0+i_1+i_2)}=H_{i_0+i_1+i_2}$, the upper ramification filtration reads:
$$H^{\nu_0} \simeq \Z /p^2 \Z \times \Z/p\Z \supset  H^{\nu_0+1}=\cdots
=H^{\nu_1} \simeq (\Z/p\Z)^2 \supset H^{\nu_1+1} =\cdots=H^{\nu_2} \simeq \Z/p\Z \supset H^{\nu_2+1} = \{0\}$$ 
with
$\nu_0=\varphi(i_0)=i_0$, $\nu_1=\varphi(i_0+i_1)=i_0+\frac{i_1}{p}$ and $\nu_2=  \varphi(i_0+i_1+i_2)= i_0+ \frac{i_1}{p} + \frac{i_2}{p^2}$.
Since $H^{p\nu_0} \supset (H^{\nu_0})^p \simeq \Z / p \Z$, we 
obtain: $H^{p\,\nu_0} \supset H^{\nu_2}$. Then, $p \, \nu_0 \leq \nu_2$, which involves  $p^2\, (p-1) \, i_0 \leq i_1 \, p +i_2$.
With such inequalities, the Hurwitz genus formula gives a new lower bound for the genus, namely: $$2 \,g=(i_0-1) (|H|-1)+i_1(|H_{i_0+1}|-1)+i_2
(|H_{i_0+i_1+1}|-1) \geq (p-1) \, (p^s \, (p^2+p+1) +(p^s+1) \, (p-1) \, p^2)$$ From
$ 2\, g \geq (p-1) \, (p^{3+s}+p^{1+s}+p^s+p^3-p^2 ) \geq (p-1) \, p^s(p^3+p)$, 
we infer the inequality: $$ \frac{|G|}{g^2}  \leq 
\frac{4}{(p^2-1)^2} \, \frac{p^{2s+3}(p+1)^2}{p^{2s}\, (p^3+p )^2} =
\frac{4}{(p^2-1)^2} \, \frac{p\, (p+1)^2}{(p^2+1 )^2}  $$ 
Since $\frac{p\, (p+1)^2}{(p^2+1 )^2}<1$ for $p \geq 2$, this contradicts condition \eqref{eq*}.
\medskip

In case iv), the lower ramification filtration of $H$ :
$$H_{i_0} \simeq \Z/p^2 \Z \times \Z/p\Z \supset  H_{i_0+1}=
\cdots=H_{i_0+i_1} \simeq (\Z/p^2\Z) \supset H_{i_0+i_1+1} =\cdots=H_{i_0+i_1+i_2} \simeq \Z / p \Z \supset \{ 0 \}. $$ 
induces the following upper ramification filtration:
$$H^{\nu_0} \simeq \Z /p^2 \Z \times \Z/p\Z \supset  H^{\nu_0+1}=\cdots
=H^{\nu_1} \simeq (\Z/p^2\Z) \supset H^{\nu_1+1} =\cdots=H^{\nu_2} \simeq \Z/p\Z \supset H^{\nu_2+1} = \{0\}.$$ 
This is almost the same situation as in case iii), except that $H_{i_0+1}$ is isomorphic to $\Z /p^2 \Z$ instead of $(\Z /p \Z )^2$. But, since the only thing that plays a part in the proof is the order of $H_{i_0+1}$ , which is the same in both cases, namely $p^2$, we conclude with the same arguments as in case iii). $\square$
\medskip

\begin{remarque}
The previous method based on the analysis of the ramification filtration of $G_2$ fails to exclude the case $G_2 \simeq \Z / p^2 \Z$ for a big action satisfying \eqref{eq*}.
Indeed, if $H:=G_2\simeq \Z / p^2 \Z$, the lower ramification filtration of $H$:
 $$H_0=\cdots=H_{i_0} \simeq \Z/p^2 \Z \supset  H_{i_0+1}=\cdots H_{i_0+i_1} \simeq \Z/p\Z \supset H_{i_0+i_1+1} = \{0 \}.$$
induces the upper ramification filtration:
$$H^0=\cdots=H^{\nu_0} \simeq \Z /p^2 \Z  \supset  H^{\nu_0+1}=\cdots
=H^{\nu_1} \simeq \Z/p\Z \supset H^{\nu_1+1} =\{0 \}.$$ 
with
$\nu_0=\varphi(i_0)=i_0$ and $\nu_1=\varphi(i_0+i_1)=i_0+\frac{i_1}{p}$.
Since $H^{p\nu_0} \supset (H^{\nu_0})^p \simeq \Z / p \Z$, we 
obtain: $p \, \nu_0  \leq \nu_1,$ hence $i_1 \geq (p-1) \, p \, i_0$.
Let $s$ be the integer defined in Remark 2.7. Then, the Hurwitz genus formula yields: $$2\, g=(i_0-1) (|H|-1)+i_1(|H_{i_0+1}|-1) \geq (p-1) \, (p^s \, (p^2+1) +p^2-p) \geq (p-1) \, p^s \, (p^2+1).$$  If we denote by $v$ the dimension of the $\F_p$-vector space $V$, we eventually get: $$ \frac{|G|}{g^2}  \leq 
 \frac{4}{(p^2-1)^2} \, \frac{p^{2+v}(p+1)^2}{p^{2s} \, (p^2+1)^2} .$$
 In this case, condition \eqref{eq*} requires $p^{1+\frac{v}{2}-s} (p+1) > p^2$. Since $\frac{v}{2} \leq s$, this implies $p+1 > p^{1+s-\frac{v}{2}} \geq p$, hence $\frac{v}{2} = s$. This means that $V= Z(Ad_f)$, where $f$ is the function defined in Remark 2.7 and $Ad_f$ its palindromic polynomial as defined in Proposition 2.5. 
Therefore, one does not obtain yet any contradiction.
 \end{remarque}

Accordingly, to exclude the cyclic cases $G_2
\simeq \Z/ p^2 \Z$ and $G_2 \simeq \Z/ p^3 \Z$ and thus complete the proof of Proposition 4.2, 
 we need to shift from a ramification point of view on $G_2$ to the embedding problem: $G_2
\subsetneq G_1$. This enables us to prove the more general result on big actions formulated in the next part.

\section{Big actions with a cyclic second ramification group $G_2$.}
 The aim of this section is to prove that there does not exist any big action whose second ramification group $G_2$ is cyclic, except for the trivial case $G_2\simeq \Z/p\Z$. 

\begin{theoreme} 
Let $(C,G)$ be a big action. If $G_2 \simeq (\Z / p^n \Z)$, then $n=1$.
\end{theoreme}

\noindent \textbf{Proof:}\\
Let $(C,G)$ be a big action with $ G_2 \simeq \Z / p^n \Z $. 
\begin{enumerate}
\item \emph{First of all, we prove that we can assume $n=2$.} \\
 Indeed, for $n >2$,  $\mathcal{H}:=G_2^{p^{n-2}}$ is a normal subgroup in $G$, strictly included in $G_2$. So Lemma 2.4.2 asserts that the pair $(C/\mathcal{H},G/\mathcal{H})$ is a big action. Besides, the second lower ramification group of $G/\mathcal{H}$ is isomorphic to $\Z/p^2\Z.$\\

\item \emph{Notation and preliminary remarks.} \\
We denote by $L:=k(C)$ the function field of $C$ and by $k(X):=L^{G_2}$ the subfield of $L$ fixed by $G_2$. Following Artin-Schreier-Witt theory (see \cite{Bo} Chap. IX, ex. 19), we define the $W_2(\F_p)$-module  $$\tilde{A}:= \frac{ \wp (W_2(L)) \cap W_2(k(X))}{ \wp (W_2(k(X)))}$$
where $W_2(L)$ denotes the ring of Witt vectors of length 2 with coordinates in $L$. 
The inclusion $k[X] \subset k(X)$ induces an injection $$ A:=\frac{ \wp (W_2(L)) \cap W_2(k[X])}{ \wp (W_2(k[X]))}\hookrightarrow  \tilde {A}$$ Since $L/L^{G_2}$ is \'etale outside $X=\infty$, it follows from \cite{Mi} (III, 4.12) that we can identify $A$ with $\tilde{A}$.
Consider the Artin-Schreier-Witt pairing: 
$$ 
\left\{
\begin{array}{ll}
G_2 \times A  \longrightarrow W_2(\F_p)\\
(g, \overline{\wp \, x}) \longrightarrow [g, \overline{\wp \, x} > := gx-x
\end{array}
\right.
$$ 
where $g \in G_2 \subset Aut_k(L)$,  $x \in L$ such that $\wp x \in k[X]$ and $\overline{\wp x}$ denotes the class of $\wp x$ mod $\wp (k[X])$. This pairing is non degenerate, which proves that, as a group, $A$ is dual to $G_2$.\\

\indent As a $\Z$-module, $A$ is generated by $(f_0(X),g_0(X))$ in $W_2(k[X])$ and then,
$L=k(X,W_0,V_0)$ with $\wp(W_0,V_0)= (f_0(X), g_0(X))$.
An exercise left to the reader shows that one can choose $f_0(X)$ and $g_0(X)$ reduced mod $\wp (k[X])$ (see definition of a reduced polynomial in section 1). We denote by $m_0$ (resp. $n_0$) the degree of $f_0$ (resp. $g_0$). Note that they are prime to $p$. The $p$-cyclic cover $L^{G_2^p}/L^{G_2}$ is parametrized by: $W_0^p-W_0=f_0(X)$. We deduce from Proposition 2.5 that $f_0(X)=XS(X)+c\,X$, where $S$ is an additive polynomial with degree $s \geq 1$ in $F$. After an homothety on $X$, we can assume $S$ to be monic. Furthermore, note that $s \geq 2$. Indeed, if $s=1$, the two inequalities established in Remark 4.4: 
$|G| \leq p^{2+2s} \leq p^4$ and $2\, g \geq (p-1) \, (p^s \, (p^2+1) +p^2-p) = (p-1) \, (p^3+p^2)$
imply $\frac{|G|}{g} \leq \frac{2\,p}{p-1} \,  \frac{p^3}{p^3+p^2} <\frac{2\,p}{p-1}$, which contradicts \eqref{N}.\\

\item \emph{The embedding problem.}\\
For any $y \in V$, the class of 
$ (f_0(X+y),g_0(X+y))$ in $A$ induces a new generating system of $A$, which means that :
\begin{equation} \label{e6}
\Z (f_0(X),g_0(X)) \,=  \, \Z (f_0(X+y),g_0(X+y)) \, \mod \wp(W_2(k[X])).
\end{equation}
As $A$ is isomorphic to $\Z /p^2 \Z$, \eqref{e6} ensures the existence of an integer $n(y)$ such that
\begin{equation} \label{e7}
(f_0(X+y),g_0(X+y))=n(y)\, (f_0(X),g_0(X)) \qquad \mod \wp(W_2(k[X]))
\end{equation}
where $n(y):=a_0(y)+b_0(y)\, p$, with $a_0(y) \in \N$, $0<a_0(y) <p$, and $b_0(y) \in \N$, $0 \leq b_0(y) <p$.
We calculate $n(y)\, (f_0(X),g_0(X))=a_0(y) \, (f_0(X),g_0(X))+  b_0(y) p \, (f_0(X),g_0(X))$.
On the one hand, $a_0(y) \, (f_0(X),g_0(X))=(a_0(y) f_0(X), a_0(y)g_0(X)+ c(a_0(y)) f_0(X))$, where $c(a_0(y))$ is given by the recursive formula:
$$
\forall \, i \in \N, \quad c(i+1)=c(i)+ \frac{1}{p} \,(1+i^p-(1+i)^p) \quad \mod \, p
$$
On the other hand, 
$$b_0(y) \,p \, (f_0(X),g_0(X))=b_0(y) \,(0,f_0(X)^p)= (0, b_0(y) f_0(X)) \, \mod \wp(W_2(k[X]))$$
As a conclusion, \eqref{e7} reads: 
\begin{equation} \label{e9}
(f_0(X+y),g_0(X+y))= (a_0(y) f_0(X), a_0(y) g_0(X)+ \ell_0(y) f_0(X)) \, \mod \wp(W_2(k[X]))
\end{equation}  where $\ell_0(y):=c(a_0(y)) +b_0(y) $. We notice that, for all $y$ in $V$, $a_0(y)=1$ mod  $p$. Indeed, the equality of the first coordinate of Witt vectors in \eqref{e9} implies 
$f_0(X+y)=a_0(y) \, f_0(X) \mod \, \wp(k[X])$. Thus, by induction, 
$f_0(X+py)=a_0(y)^p \, f_0(X) \mod \, \wp(k[X])$. Since $V$ is an elementary abelian $p$-group, $f_0(X+py)=f_0(X)$, which involves: $a_0(y)^p=1$ mod $p$ and $a_0(y)=1$ mod $p$.
So, \eqref{e9} becomes:
\begin{equation} \label{e10}
 (f_0(X+y),g_0(X+y))= (f_0(X),g_0(X)+ \ell_0(y) f_0(X)) +(P^p(X),Q^p(X))-(P(X),Q(X))
\end{equation}
with $P(X)$ and $Q(X)$ polynomials of $k[X]$. In order to circumvent the problem related to the special formula giving the opposite of Witt vectors for $p=2$, we would rather write \eqref{e10} as follows:
\begin{equation} \label{e11}
(f_0(X+y),g_0(X+y))+(P(X),Q(X))= (f_0(X),g_0(X)+ \ell_0(y) \,f_0(X)) +(P(X)^p,Q(X)^p) 
\end{equation}
The first coordinate of \eqref{e11} reads:
\begin{equation} \label{e12}
f_0(X+y)+P(X)=f_0(X)+P(X)^p 
\end{equation}
On the second coordinate of \eqref{e11}, the addition law in the ring of Witt vectors gives the following equality in $k[X]$:
\begin{equation} \label{e13}
g_0(X+y)+Q(X)+\psi(f_0(X+y),P(X))=g_0(X)+\ell_0(y)\, f_0(X)+Q(X)^p+ \psi(f_0(X),P(X)^p)
\end{equation}
where $\psi$ is defined as follows: $$\psi(a,b):= \frac{1}{p}\, (a^p+b^p-(a+b)^p)
= \frac{-1}{p}\, \sum_{i=1}^{p-1}\, \binom pi \, a^i\,  b^{p-i}= \sum_{i=1}^{p-1}\, \frac{(-1)^i}{i} \, a^i \,  b^{p-i}\quad  \mod \, p$$
As a consequence, \eqref{e13} gives:
\begin{equation} \label{e14}
\Delta_y (g_0):= g_0(X+y)-g_0(X)=\ell_0(y) \, f_0(X)+ \delta \qquad \mod \, \wp (k[X]) 
\end{equation}
with 
$$\begin{array}{ll}
\delta&:=\psi(f_0(X),P(X)^p)-\psi(f_0(X+y),P(X))\\
&= \sum_{i=1}^{p-1} \,\frac{(-1)^i}{i}\, \{ f_0(X)^i \, P(X)^{p(p-i)}-f_0(X+y)^i \,P(X)^{p-i} \}
\end{array}$$

\begin{lemme}
 With the notation defined above, $\delta$ is equal to:
 \begin{equation} \label{e15}
 \delta =\sum_{i=1}^{p-1} \frac{(-1)^i}{i} \, y^{p-i} X^{i+p^{s+1}} +\mbox {lower  degree terms in  X }
\end{equation}
\end{lemme}

\noindent \textbf{Proof:}
We search for the monomials in $\delta$ that have degree in $X$ greater or equal to $p^{s+1}+1$. We first focus on $f_0(X)^i \, P(X)^{p(p-i)}$.
We can infer from equality \eqref{e12} that $P(X)$ has degree $p^{s-1}$ and that its leading coefficient is $y^{1/p}$. Furthermore, \cite{LM} (proof of Prop. 8-1) shows that $P(X)-P(0)$ is an additive polynomial. So, we can write: $P(X)=y^{1/p} \, X ^{p^{s-1}} +P_1(X)$,
where $P_1(X)$ is a polynomial of $k[X]$ with degree at most $p^{s-2}$.
Then, for all $i$ in $\{1,\cdots,p-1\}$, $f_0(X)^i\, P(X)^{p\, (p-i)}= f_0(X)^i \, (y\, X^{p^s}+P_1(X)^p)^{p-i}= f_0(X)^i \, (\sum_{j=0}^{p-i} \binom {p-i}j \, y^j \, X^{jp^s} \, P_1(X)^{p(p-i-j)})$. Since $f_0(X)$ has degree: $1+p^s$, this gives in $\delta$ a monomial of degree at most:
$i\, (1+p^s)+j\, p^s +p\, (p-i-j)\, p^{s-2}=p^s+(i+j)\, (p-1) \,p^{s-1}+i.$
If $j \leq p-i-1$, this degree is at most: 
$p^s+(p-1)^2\,p^{s-1}+i=(p-1)\,p^s+p^{s-1}+i$, which is strictly lower than $p^{s+1}+1$, for $s \geq 2$ and $1 \leq i \leq p-1$ . 
As a consequence, the monomials of degree greater or equal to $p^{s+1}+1$ can only occur when the index $j$ is equal to $p-i$, namely in $f_0(X)^i \,  y^{p-i} \, X^{p^s(p-i)}$. As $f_0(X)=X \,S(X)+c \, X$, where $S$ is a monic additive polynomial of degree $s$ in $F$, $f_0$ reads: $f_0(X)=X^{1+p^s} +P_2(X)$ where $P_2(X)$ is a polynomial in $k[X]$ with degree at most $1+p^{s-1}$. 
Then, for all $i$ in $\{1,\cdots,p-1\}$, $f_0(X)^i \,  y^{p-i} \, X^{p^s(p-i)}=  y^{p-i} \, X^{p^s(p-i)} \, (\sum_{k=0}^i \binom ik X^{(1+p^s)j}\, P_2(X)^{i-k})$.
Accordingly, we get a monomial of degree at most: $ p^s\, (p-i) +k\, (1+p^s)+(i-k) \,(1+p^{s-1})=p^s \, (p-i) +i \, (1+p^{s-1} ) +k \, (p^s-p^{s-1})$. When $0 \leq k \leq i-1$, the maximal degree obtained in this way is $i +  p^{s-1}-p^s+p^{s+1}$ which is stricly lower than $p^{s+1}+1$. Therefore, for all $i$ in $\{1,\cdots,p-1\}$, the only contibution to take into account is $k=i$, which produces in $\delta$ the sum: $ \sum_{i=1}^{p-1} \frac{(-1)^i}{i}\, y^{p-i} X^{i+p^{s+1}}$.\\
\indent We now search for monomials with degree greater or equal to $p^{s+1}+1$ in the second part of $\delta$, namely: $f_0(X+y)^i \,P(X)^{p-i}$. This has degree at most: $i\, (1+p^s) +(p-i)\,p^{s-1}=i\,p^s+(p-i)\, p^{s-1}+i$, which is strictly lower than $p^{s+1}+1$, for $s\geq 2$ and $1 \leq i \leq p-1$. Therefore, $f_0(X+y)^i \,P(X)^{p-i}$ does not give any monomial in $\delta$ with degree greater or equal to $p^{s+1}+1$. Thus, we get the expected formula. $\square$ $\square$ 
\medskip

\item \emph{We notice that $g_0(X)$ cannot be of the form $X\, \Sigma(X)+\gamma\,X$, with $\Sigma \in k\{F\}$ and $\gamma \in k$.}\\ Otherwise, the left-hand side of \eqref{e14} reads: $\Delta_y(g_0):=
g_0(X+y)-g_0(X)=X\,\Sigma(y)+y\,\Sigma(X)+y\,\Sigma(y)+\gamma \,y$, which only gives a linear contribution in $X$ after reduction mod $\wp(k[X])$. By Lemma 5.2, $deg \, f_0=1+p^s < deg\, \delta= p^{s+1}+p-1$, which involves that the degree of the right-hand side of \eqref{e14} is $p-1+p^{s+1} >1$, hence a contradiction.\\
\indent Therefore, we can define an integer $a \leq n_0=deg \, g_0$ such that $X^{a}$ is the monomial of $g_0(X)$ with highest degree which is not of the form $1+p^n$, with $n\in \N$. Note that since $g_0$ is reduced mod $\wp(k[X])$, $a \not \equiv 0$ mod $p$.
We also notice that the monomials in $g_0(X)$ with degree strictly greater than $a$ are of the form $X^{1+p^n}$, and so, as explained above, they only give linear monomials in $\Delta_y(g_0)$ mod  $\wp(k[X])$. Therefore, after reduction mod $\wp(k[X])$, the degree of the left-hand side of \eqref{e14} is at most $a-1$. Since the degree of the right-hand side is $p^{s+1}+p-1$, it follows that: 
\begin{equation} \label{a}
a-1 \geq p^{s+1}+p-1 \qquad
\end{equation}

\item \emph{We show that $p$ divides $a-1$.}\\
Assume that $p$ does not divide $a-1$. In this case, the monomial $X^{a-1}$ is reduced mod $\wp(k[X])$. Since the monomials of $g_0(X)$ with degree strictly greater than $a$ only give a linear contribution in $\Delta_y(g_0)$ mod $\wp \, (k[X])$, \eqref{e14} reads as follows, for all $y$ in $V$:
$$c_a(g_0) \, a\, y X^{a-1}+ \mbox{lower degree terms } = -y \, X^{p^{s+1}+p-1} + 
\mbox{ lower degree  terms} \mod \, \wp\, (k[X])
$$ where $c_a(g_0)\neq 0$ denotes the coefficient of $X^a$ in $g_0$.
 If $a-1 > p^{s+1}+p-1$, the coefficient $c_a(g_0)\, a\, y=0$, for all $y$ in V. Since $a
\neq 0$ mod $p$, it leads to $V=\{0\}$, so $G_1=G_2$, which is impossible for a big action
(see Proposition 2.2.1). We gather from \eqref{a} that $a-1= p^{s+1}+p-1$, which contradicts: $a \neq 0$ mod $p$.\\
\indent Thus, $p$ divides $a-1$. So, we can write $a=1+\lambda\, p^t$, with $t>0$, $\lambda$ prime to $p$ and $\lambda \geq 2$ because of the definition of $a$. We also define $j_0:=a-p^t=1+(\lambda-1)\, p^t$. Note that $pj_0 > a$. Indeed, 
$$pj_0 \leq a \Leftrightarrow p(1+(\lambda-1)p^t) \leq 1+ \lambda\,p^t \Leftrightarrow \lambda \leq \frac{1-p+p^{t+1}}{p^t(p-1)}=\frac{-1}{p^t}+ \frac{p}{p-1} < \frac{p}{p-1} \leq 2$$ which is impossible since $\lambda  \geq 2$.
 
\item \emph{We determine the coefficient of $X^{j_0}$ in the left hand-side of \eqref{e14}.}\\ 
Since $p$ does not divide $j_0$, the monomial $X^{j_0}$ is reduced mod $\wp(k[X])$.
 In the left-hand side of \eqref{e14}, namely $\Delta_y(g_0)$ mod $\wp(k[X])$, the monomial 
$X^{j_0}$ comes from monomials of $g_0(X)$ of the form: $X^b$, with $b$ in $\{j_0+1,\cdots,a\}$. As a matter of fact, the monomials of $g_0(X)$ with degree strictly greater than $a$ only give a linear contribution mod $\wp(k[X])$, whereas $j_0=1+(\lambda-1)\, p^t>1$. For all $b \in \{j_0+1,\cdots, a\}$, the monomial $X^b$ of $g_0(X)$ generates $\binom b{j_0}\,  y^{b-j_0} \, X^{j_0}$ in $\Delta_y(g_0)$. Since $p\, j_0 > a \geq b$ (see above), these monomials $X^b$ do not produce any $X^{j_0\, p^n}$, with $n \geq 1$, which would also give $X^{j_0}$ after reduction mod $\wp(k[X])$. It follows that the coefficient of $X^{j_0}$ in the left-hand side of \eqref{e14} is $T(y)$ with $T(Y):= \sum _{b=j_0+1}^a \, c_b(g_0)\, \binom b{j_0}\, Y^{b-j_0}$, where $c_b(g_0)$ denotes the coefficient of $X^b$ in $g_0(X)$. As the coefficient of $Y^{a-j_0}$ in $T(Y)$ is $c_a(g_0) \, \binom a{j_0}=c_a(g_0)\, \binom {1+\lambda p^t}{1+(\lambda-1)p^t} \equiv c_a(g_0) \, \lambda \not \equiv 0 \, \mod \, p$, the polynomial $T(Y)$ has degree $a-j_0=p^t$. \\

\item \emph{We identify with the coefficient of $X^{j_0}$ in the right-hand side of \eqref{e14} and gather a contradiction.}\\
We first assume that the monomial $X^{j_0}$ does not occur in the right-hand side of \eqref{e14}. Then, $T(y)=0$ for all $y$ in $V$, which means that $V$ is included in the set of roots of $T$. Thus, $|V| \leq p^t$.
To compute the genus $g$, put $M_0:=m_0$ and $M_1:=\max \{ p \,m_0, \,n _0 \}$. Then, by \cite{Ga}, the Hurwitz genus formula applied to $C \rightarrow C/G_2 \simeq \p_k^1$  yields: $2\, (g-1)=2\, |G_2|\, (g_{C/G_2}-1) +d=-2\,p^2+d$, with $d:=(p-1)\,(M_0+1)+p\,(p-1) \,(M_1+1)$. From $p \, m_0= p\, (p^s+1)=p^{s+1}+p$ and $p^{s+1}+p-1 < n_0$, we infer $M_1=n_0$. Moreover, since $n_0 \geq a =1+\lambda \, p^t \geq 1+2\,p^t >2\,p^t$,  we obtain a lower bound for the genus:  $2\, g=(p-1)\,p\, (n_0-1+p^{s-1}) \geq 2\,p^{t+1} \, (p-1)$.
As $|G|=|G_2||V|\leq p^{2+t}$, it entails an inequality $\frac{|G|}{g} \leq \frac{2\:p}{p-1} \, \frac{p^{1+t}}{2\,p^{1+t}}= \frac{1}{2} \, \frac{2\,p}{p-1}$
which contradicts \eqref{N}.\\
\indent As a consequence, the monomial $X^{j_0}$ appears in the right-hand side of \eqref{e14}, which implies that $j_0 \leq p^{s+1}+p-1$. Using \eqref{a}, we get:
$ j_0=1+(\lambda-1)\,p^t \leq p^{s+1}+p-1 < a=1+\lambda \,p^t $. 
 This yields:
\begin{equation} \label{t2}
 \lambda-1 \leq p^{s+1-t} + \frac{p-2}{p^t} < \lambda
\end{equation}
 If $ s+1-t \leq -1$, since $t \geq 1$, \eqref{t2} gives: $\lambda-1 \leq \frac{1}{p} +\frac{p-2}{p} <1$, which contradicts $\lambda \geq 2$. It follows that $s+1-t \geq 0$. Then, \eqref{t2} combined with: $0 \leq \frac{p-2}{p^t} <1$  leads to $\lambda-1=p^{s+1-t}$. We gather that $j_0=1+(\lambda-1)\, p^t =1+p^{s+1} > deg \, f_0=1+p^s$. Therefore, in the right-hand side of \eqref{e14},  the monomial $X^{j_0}=X^{1+p^{s+1}}$ only occurs in $\delta$. By Lemma 5.2, the coefficient of $X^{j_0}=X^{1+p^{s+1}}$ in $\delta$ is $-y^{p-1}$. By equating the coefficient of $X^{j_0}$ in each side of \eqref{e14}, we get
$T(y)=-y^{p-1}$, for all $y$ in $V$. Put $\tilde{T}(Y):=
T(Y)+Y^{p-1}$. Since $deg \, T =p^t > p-1$, the polynomial $\tilde{T}$ has still degree $p^t$ and satisfies: $\tilde{T}(y)=0$ for all $y$ in $V$. Once again, it leads to $|V| \leq p^t$, which contradicts \eqref{N} as above. $\square$
\end{enumerate} 

\medskip 

Therefore, when $(C,G)$ is a big action, $G_2 \simeq (\Z/p^n \Z)$  implies $n=1$.
More generally, if $G_2$ is abelian of exponent $p^n$, with $n \geq 2$, there exists an index $p$-subgroup of $G_2^p$, say $H$, normal in $G$ such that the pair $(C/H, G/H)$ is a big action with $(G/H)_2=G_2/H \simeq \Z/p^2 \Z \times (\Z/p\Z)^t$, with $t \in \N^*$. A natural question is to search for a lower bound on the $p$-rank: $t$ depending on the genus $g$ of the curve.  As seen in the proof of Theorem 5.1, the difficulty lies in the embedding problem, i.e. in finding an extension which is stable under the translations by $V$. In the next section, we exhibit big actions with $G_2$ abelian of exponent at least $p^2$. In particular, we construct big actions $(C,G)$ with $G_2 \simeq \Z/p^2\Z \times (\Z/p\Z)^t$ where $t= O(log_p \, g)$.

\section{Examples of big actions with $G_2$ abelian of exponent strictly greater than $p$.}
 In characteristic $0$, an anologue of big actions is given by the actions of a finite group $G$ on a compact Riemann surface $C$ with genus $g_C \geq 2$ such that $|G| = 84(g_C-1)$. Such a curve $C$ is called a \emph{Hurwitz curve} and such a group $G$ a \emph{Hurwitz group} (cf. \cite{Co}). In particular, the lowest genus Hurwitz curves are the Klein's quartic with $G \simeq PSL_2(\F_7)$ (cf. \cite{El2}) and the Fricke-Macbeath curve with genus $7$ and $G \simeq PSL_2(\F_8)$ (cf. \cite{Mc65}). \\
\indent Let $C$ be a Hurwitz curve with genus $g_c$. Let $n \geq 2$ be an integer and let $C_n$ be the maximal unramified Galois cover whose group is abelian, with exponent $n$. The Galois group of the cover $C_n/ C$ is isomorphic to $(\Z / n \Z)^{2g_C}.$
We infer from the unicity of $C_n$ that the $\C$-automorphims of $C$ have $n^{2g_c}$ prolongations to $C_n$. Therefore, $g_{C_n}-1
= n^{2g}(g_{C}-1)$. Consequently, $C_n$ is still a Hurwitz curve (see \cite{Mc}).\\
\indent Now, let $(C,G)$ be a big action. Then $C \rightarrow C/G$ is an \'etale cover of the affine line whose group is a $p$-group.
From the Deuring-Shafarevich formula (see e.g. \cite{Bou}), it follows that the Hasse-Witt invariant of $C$ is zero. This means that there are no nontrivial connected \'etale Galois covers of $C$ with group a $p$-group. Therefore, if we want to generalize the method mentionned above to produce Galois covers of $C$ corresponding to big actions, it is necessary to introduce ramification.
A means to do so is to consider ray class fields of function fields, as studied by K. Lauter \cite{Lau} and R. Auer \cite{Au}. Since the cover $C \rightarrow C/G_2$ is an \'etale cover of the affine line $Spec \, k[X]$ totally ramified at $\infty$, we focus on the special case of ray class fields of the rational function field $\F_q(X)$, where $q=p^e$ (see \cite{Au}, III.8).
Such ray class fields allows us to produce families of big actions $(C,G)$ (where $C$ is defined over $k=\F_p^{alg}$) with specific conditions imposed on ramification and endowed with an abelian $G_2$ of exponent as large as we want.

\begin{defi} (\cite{Au}, Part II)
Let $K$ be the rational function field: $\F_q(X)$, with $q=p^e$ and $e \in \N^*$.
Let $S$ be the set of all finite rational places, namely $\{ (X-y),\, y \in \F_q \}$.
Let $m \geq 0$ be an integer.
Fix  $K^{alg}$ an algebraic closure of $K$ in which all extensions of $K$ are assumed to lie.
We  define $K_S^m \subset  K^{alg}$ as the largest abelian extension $L / K$ with  conductor $\leq   m \,\infty$, such that every  place in $S$ splits completely in $L$.
\end{defi}

\begin{remarque}
\begin{enumerate}
\item  We define the splitting set of any finite Galois extension $L/K$, denoted by $S(L)$, as the set consisting of the places of $K$ that split completely in $L$. If $K_S^m/K$ is the extension defined in Definition 6.1, then $S\subset S(K_S^m)$.
\item In what follows, we only consider finite Galois extensions $L/K$ that are unramified outside $X=\infty$ and (totally) ramified at $X=\infty$. Therefore, the support of the conductor of $L/K$ is  reduced to the place $\infty$. So, we systematically confuse the conductor $m\,\infty$ with its degree $m$.
\item We could more generally define $K_S^m$ for $S$ a non-empty subset of the finite rational places, i.e. $S:= \{ (X-y), \, y \in V \subset \F_q\}$. However, to get big actions, it is necessary to consider the case where $V$ is a subgroup of $\F_q$. In what follows, we focus on the case $V=\F_q$, as announced in Definition 6.1.
\end{enumerate}
\end{remarque}

\begin{remarque}  We keep the notation of Definition 6.1. 
\begin{enumerate}
\item The existence of the extension $K_S^m / K$ is based on global class field theory (see \cite{Au}, Part II). 
\item $K_S^m/K$ is a finite abelian extension whose full constant field is $\F_q$. 
\item The reason for Lauter and Auer's interest in such ray class fields is that they provide for examples of global function fields with many rational places, or what amounts to the same, of algebraic curves with many rational points. Indeed, let $C(m) / \F_q$ be the nonsingular projective curve with function field $K_S^m$. If we denote by $N_m:= |C(m)(\F_q)|$ the number of $\F_q$-rational points on the curve $C(m)$, then: $ N_m=1+q \,[K_S^m:K]$.
The main difficulty lies in computing $[K_S^m:K]$. We first wonder when $K_S^m$ coincide with $K$. Here are partial answers.
\item Let $q=p^e$, with $ e\in \N$. If $e$ is even, put $r:=\sqrt{q}$ and if $e$ is odd, put $r:=\sqrt{qp}$.
Then, for all $i$ in $\{0,\cdots, r+1\}$, $K_S^i=K=\F_q(X)$. (see \cite{Au}, III, Lemma 8.7 and formula (13)). Note that the previous estimate $ N_m=1+q \,[K_S^m:K]$, combined with the Hasse-Weil bound (see e.g. \cite{St93} V.2.3), furnishes another proof of $K_S^i=K$ when $i <1+r$.
\item More generally, Lauter displays a method to compute the degree of the extension $K_S^m/K$ via a formula giving the order of its Galois group: $G_S(m)$ (see \cite{Lau}, Thm. 1). Her proof consists in starting from the following presentation of $G_S(m)$:
$$ G_S(m) \simeq \frac{1+Z\, \F_q [[Z]]} {< 1+Z^m\, \F_q[[Z]], 1-yZ,\, y \in \F_q >}$$
where $Z=X^{-1}$, which indicates that $G_S(m)$ is an abelian finite $p$-group. Then, she transforms the multiplicative structure of the group into an additive group of generalized Witt vectors.
In particular, she deduces from this theorem  the smallest conductor $m$ such that $G_S(m)$ has exponent stricly greater than $p$ (see next proposition).
\end{enumerate}
\end{remarque}

\begin{proposition} (\cite{Lau}, Prop. 4)
We keep the notation defined above.
If $q=p^e$, the smallest conductor $m$ for which the group $G_S(m)$ is not of exponent $p$ is $m_2:= p^{\lceil e/2\rceil +1} +p+1$, where $\lceil e/2 \rceil$ denotes the upper integer part of $e/2$. 
\end{proposition}

We now emphasize the link with big actions. Let $F$ be a function field with full constant field $\F_q$. Let $C/\F_q$ be the smooth projective curve whose function field is $F$ and $C^{alg}:= C\times_{\F_q} k$ with $k=\F_p^{alg}$. If $G$ is a finite $p$-subgroup of $Aut_{\F_q} C$, then $G$ can be identified with a subgroup of $Aut_{k} C^{alg}$.
In this case, $(C^{alg}, G)$ is a big action if and only if $g_{C^{alg}}=g_C >0$ and $\frac{|G|}{g_C} >\frac{2\,p}{p-1}$.
For convenience, in the sequel, we shall say that $(C,G)$ is a big action if $(C^{alg}, G)$ is a big action.\\
\indent In what follows, we consider the curve $C(m)/\F_q$ whose function field is $K_S^m$ and, starting from this, we construct a $p$-group $G(m)$ acting on $C(m)$ by extending the translations $X \rightarrow X+y$, with $y \in \F_q$.
In particular, we obtain an upper bound for the genus of $C(m)$, which allows us to circumvent the problem related to the computation of the degree $[K_S^m:K]$ when checking whether $(C(m), G(m))$ is a big action.

\begin{proposition}
We keep the notation defined above.
\begin{enumerate}
\item Let $C(m) / \F_q$ be the nonsingular projective curve with function field $K_S^m$.
Then, the group of translations:  $X \rightarrow X+y$, $y\in \F_q$, extends to a $p$-group of $\F_q$-automorphisms of $C(m)$, say $G(m)$, 
with the following exact sequence:
$$ 0 \longrightarrow G_S(m) \longrightarrow G(m) \longrightarrow \F_q \longrightarrow 0$$
\item Let $L$ be an intermediate field of $K_S^m/K$. Assume $L=(K_S^m)^H$, i.e. the extension $L/K$ is Galois with group: $G_S(m)/H$. For all $i \geq 0$, we define $L^i$ as the $i$-th upper ramification field of $L$, i.e. the subfield of $L$ fixed by the $i$-th upper ramification group of $G_S(m)/H$ at $\infty$: $G_S^i(m)H/H$, where $G_S^i(m)$ denotes the $i$-th upper ramification group of $G_S(m)$ at $\infty$. Then, 
$$ \forall \, i \geq 0, \quad L^i= L \cap K_S^{i}$$
In particular, when $L=K_S^m$ and $ i \leq m$, 
$L^i=K_S^i$, i.e. $G_S^i(m)=Gal(K_S^m/K_S^i)$. 
\item Let $L$ be an intermediate field of $K_S^m/K$. Define $n:=\min \{i \in \N, L \subset K_S^i\}$. Then, the genus of the extension $L/K$ is given by the formula: 
$$g_{L} = 1+ [L:K] \,(-1 + \frac{n}{2}) -\frac{1}{2}\, \sum_{j=0}^{n-1} \,[L \cap K_S^j:K]$$
where the sum is empty for $n=0$.\\
In particular, $g_L=0$ if and only if $n:=\min \{i \in \N, L \subset K_S^i\}=0$.\\
Note that if $n >0$, then $g_L < [L:K] \, (-1+\frac{n}{2})$. 
\item If $m \geq r+2$,  $\frac{|G(m)|}{g_{K_S^{m}}} >\frac{q}{-1+\frac{m}{2}}$. It follows that if $\frac{q}{-1+\frac{m}{2}} \geq \frac{2\,p}{p-1}$, the pair $(C(m), G(m))$ is a big action. In this case, the second lower ramification group of $G(m)$: $G_2(m)$, is equal to $G_S(m)$. 
In particular, for $p >2$, (resp. $p=2$), if $e\geq 4$ (resp. $e \geq 6$) and if $m_2$ is the integer defined in Proposition 6.4, the pair $(C(m_2), G(m_2))$ is a big action whose second  ramification group: $G_S(m_2)$, is abelian of exponent $p^2$.
\end{enumerate}
\end{proposition}

\noindent \textbf{Proof:}
\begin{enumerate}
\item The set $S$ is globally invariant under the translations: $X \rightarrow X+y$, $y\in \F_q$. That is the same for $\infty$, so the translations by $\F_q$ do not change the conditions imposed on ramification. 
As a consequence, owing to the maximality and the unicity of $K_S^m$, they can be extended to $\F_q$-automorphisms of $K_S^m$. This proves the first assertion. 
\item The second point directly derives from \cite{Au} (II, Thm. 5.8). 
\item The genus formula is obtained by combining the preceding results, the Hurwitz genus formula and the Discriminant formula (see \cite{Au}, I, 3.7). 
Now assume that $n=0$. Then, $L\subset K_S^0=\F_q(X)$ and $g_L=0$.
Conversely, assume $g_L=0$. If $n \neq 0$, Remark 6.3.4 implies that $n \geq r+2 \geq 3$. Using the preceding formula and Remark 6.3.4, $g_L=0$ reads:
$$ 2+ (n-2)\, [L:K] = \sum_{j=0}^{n-1} \, [K_S^j \cap L:K]= 2+ \sum_{j=2}^{n-1} \,  [K_S^j \cap L:K] \leq 2+ (n-2)\, [L:K]$$
It follows that, for all $j$ in $\{2,\cdots,n-1\}$, $K_S^j \cap L=L$. In particular, $L \subset K_S^2=K_S^0$, hence a contradiction. Finally, since $n>0$ implies $n \geq 3$ and since $K=K_S^0=K_S^1$, one notices that
$$g_L=[L:K] \,(-1 + \frac{n}{2}) -\frac{1}{2}\, \sum_{j=2}^{n-1} \,[L \cap K_S^j:K] <[L:K]\, (-1+\frac{n}{2}) $$
\item Assume that $m \geq r+2$. We gather from Remark 6.3.4 that $n:=\min \{i \in \N, K_S^m \subset K_S^i\} \geq r+2 \geq 3$. Then, it follows from the previous point that
$$g_{K_S^m} <[K_S^m:K]\, (-1+\frac{n}{2}) \leq [K_S^m:K] \,(-1 + \frac{m}{2}) $$
As $|G(m)|=q [K_S^m:K]$, we deduce the expected inequality. In particular, when $\frac{q}{-1+\frac{m}{2}} > \frac{2\,p}{p-1}$, the pair $(C(m), G(m))$ is a big action.
It remains to show that, in this case, $G_2(m)$ is equal to $G_S(m)$. Lemma 2.4.2 first proves that $G_S(m) \supset G_2(m)$. Let $L:=(K_S^m)^{G_2(m)}$ be the subfield of $L$ fixed by $G_2(m)$ and define $n:=\min \{i \in \N, L \subset K_S^i\}$. Assume $G_S(m) \supsetneq G_2(m)$. Then $L \supsetneq (K_S^m)^{G_S(m)}=K$. We infer from Remark 6.3.4 that $n \geq r+2$, which proves, using the previous point, that $g_L>0$. But, since $(C(m), G(m))$ is a big action, $C/G_2(m) \simeq \p_k^1$, so $g_L=0$, hence a contradiction. We eventually explain the last statement.
By Proposition 6.5.2, $G_S^{m_2-1}(m_2)=Gal(K_S^{m_2}/K_S^{m_2-1})$, which induces the following exact sequence:
$$ 0 \longrightarrow G_S^{m_2-1}(m_2) \longrightarrow G_S(m_2) \longrightarrow G_S(m_2-1) \longrightarrow 0$$
We infer from Proposition 6.4 that $G_S(m_2-1)$ has exponent $p$ whereas the exponent of $G_S(m_2)$ is at least $p^2$. It follows that $G_S^{m_2-1}(m_2)$ cannot be trivial.
Since $G_S^{m_2}(m_2)=\{0\}$ (use Proposition 6.5.2), we deduce from the elementary properties of the ramification groups that $G_S^{m_2-1}(m_2)$ is $p$-elementary abelian. Therefore, $G_S(m_2)$ has exponent smaller than $p^2$ and the claim follows. $\square$
\end{enumerate}

\begin{remarque}
Let $N_m$ be the number of $\F_q$-rational points on the curve $C(m)$ as defined in Remark 6.3.3. Then, $N_m=1+q \, |G_S(m)|=1+ |G(m)|$. This highlights the equivalence of the two ratios: $\frac{|G(m)|}{g_{C(m)}}$ and $\frac{N_m}{g_{C(m)}}$. In particular, this equivalence emphasizes the link between the problem of big actions and the search of algebraic curves with many rational points.
\end{remarque}

As seen in Remark 6.3.4, $K_S^i=K$ for all $i$ in $\{0,\cdots,r+1\}$, where $r=\sqrt{q}$ or $\sqrt{qp}$ according to whether $q$ is a square or not. 
The following extensions $K_S^m$, for $m \geq r+2$, are partially parametrized, at least for the first ones, in \cite{Au} (Prop. 8.9). In the table below, we exhibit a complete description of the extensions $K_S^m$ for $m$ varying from $0$ to $m_2=p^{\lceil e/2\rceil +1}+p+1$, in the special case $p=5$ and $e=4$. This involves $q=p^e=625$, $s=e/2=2$, $r=p^s=25$ and $m_2=131$. The table below should suggest the general method to parametrize such extensions.

\bigskip

\begin{tabular}{|c|c|c|}
\hline
 conductor $m$ & $[K_S^m:K]$ & New equations\\
\hline
$0 \leq m \leq r+1=26$ & $1$ & \\
\hline
$r+2=27 \leq m \leq 2r+1=51$ & $5^2$ & $W_{0}^r+W_{0}=X^{1+r}$\\
\hline
$m=2r+2=52$ & $5^6$ & $W_{1}^q-W_{1}=X^{2r}\, (X^q-X)$\\
\hline
 $2r+3=53 \leq m \leq 3r+1=76$& $5^8$ & $W_{2}^r+W_{2}=X^{2(1+r)}$\\
\hline
$m=3r+2=77$ & $5^{12}$ & $W_{3}^q-W_{3}=X^{3r}\, (X^q-X)$\\
\hline
$m=3r+3=78$ & $5^{16}$ & $W_{4}^q-W_{4}=X^{3r}\, (X^{2q}-X^2)$\\
\hline
 $3r+4=79 \leq m \leq 4r+1=101$& $5^{18}$ & $W_{5}^r+W_{5}=X^{3(1+r)}$\\
\hline
$m=4r+2=102$ & $5^{22}$ & $W_{6}^q-W_{6}=X^{4r}\, (X^q-X)$\\
\hline
$m=4r+3=103$ & $5^{26}$ & $W_{7}^q-W_{7}=X^{4r}\, (X^{2q}-X^2)$\\
\hline
$m=4r+4=104$ & $5^{30}$ & $W_{8}^q-W_{8}=X^{4r}\, (X^{3q}-X^3)$\\
\hline
 $4r+5=105 \leq m \leq 5r+1=126$& $5^{32}$ & $W_{9}^r+W_{9}=X^{4(1+r)}$\\
\hline
$m=5r+2=127$ & $5^{36}$ & $W_{10}^q-W_{10}=X^{5r}\, (X^{q}-X)$\\
\hline
$m=5r+3=128$ & $5^{40}$ & $W_{11}^q-W_{11}=X^{5r}\, (X^{2q}-X^2)$\\
\hline
$m=5r+4=129$ & $5^{44}$ & $W_{12}^q-W_{12}=X^{5r}\, (X^{3q}-X^3)$\\
\hline
$m=5r+5=130$ & $5^{48}$ & $W_{13}^q-W_{13}=X^{5r}\, (X^{4q}-X^4)$\\
\hline
$m=m_2=131$ & $5^{50}$ & $[W_{0}, W_{14}]^r+[W_{0}, W_{14}]=[X^{1+r},0] $\\
\hline
\end{tabular}

\bigskip

In this case, 
\begin{equation} \label{quotient}
\frac{|G(m_2)|}{g_{K_S^{m_2}}} \simeq 9,6929 \cdots
\end{equation}

\medskip

\noindent \textbf{Comments on the construction of the table:}
For all $i$ in $\{0,\cdots,14\}$, put $L_i:=K(W_0, \cdots, W_i)$.
\begin{enumerate}
\item We first prove that the splitting set of each extension $K(W_i)/K$ (see Remark 6.2.1) contains $S$. Indeed, fix $y$ in $\F_q$ and call $P_y$ the corresponding place in $S$: $(X-y)$. We have to distinguish three cases.
By \cite{St93} (Prop. VI. 4.1), $P_y$ completely splits in the extension $K(W)/K$, where $W^r+W=X^{u\,(1+r)}$, with $1 \leq u \leq 4$, if the polynomial $T^r+T-y^{u\,(1+r)}$ has a root in $K$, which is true since $y^{u(1+r)}=(F^s+I) \, (\frac{1}{2}\, y^{u(1+r)})$.
Likewise, $P_y$ completely splits in the extension $K(W)/K$, where $W^q-W=X^{u\,r} \, (X^{v\,q}-X^v)$, with $1 \leq v < u \leq 5$, since $y^{vq}-y^v=0$.
Finally, $P_y$ completely splits in the extension $K(W,\tilde{W})/K$, where $[W,\tilde{W}]^r+[W,\tilde{W}]=[X^{1+r},0]$, since $[y^{1+r},0]= (F^s+I)\, [\frac{1}{2} \,y^{1+r},-\frac{2^p-2}{4p}\, y^{(1+r)\,p}]$.
To conclude, we remark that $L_i=L_{i-1} \, K(W_i)$ for all $i$ in $\{1,\cdots,14\}$. Then, $S(L_i)=S(L_{i-1}) \cap S(K(W_i))$ (cf.  \cite{Au}, Cor. 3.2.b),  which allows us to gather, by induction on $i$, that the splitting set of each $L_i$ contains $S$.
\item We now compute the conductor $m(K(W_i))$ of each extension $K(W_i)/K$. As above, we have to distinguish three kinds of extensions. First, the extension $K(W)/K$, where $W^r+W=X^{u\,(1+r)}$, with $1 \leq u \leq 4$, has conductor $ur+u+1$
 (see \cite{Au}, Prop. 8.9.a). Besides, the extension $K(W)/K$, where $W^q-W=X^{u\,r} \, (X^{v\,q}-X^v)$, with $1 \leq v < u \leq 5$, has conductor $ur+v+1$ (see \cite{Au}, Prop. 8.9.b). Finally, the conductor of the extension $K(W,\tilde{W})/K$, where $[W,\tilde{W}]^r+[W,\tilde{W}]=[X^{1+r},0]$ is given by the formula:
$1+ \max \{p(1+r),0\} = 1+p+p^{s+1}=m_2$ (see \cite{Ga}, Thm. 1.1).
As a conclusion, since $m(L_i)= \max \{m(L_{i-1}), m(K(W_i)) \}$ (cf. \cite{Au}, Cor. 3.2.b), an induction on $i$ allows us to obtain the expected conductor for $L_i$.
\item We gather from the two first points the inclusions:
$K(W_0) \subset K_S^{27}$,  $K(W_0,W_1) \subset K_S^{52},  \cdots$ \\$  K(W_0,\cdots, W_{14}) \subset K_S^{m_2}$.
Equality is eventually obtained by calculating the degree of each extension $K_S^m/K$ via \cite{Lau} (Thm. 1) or \cite{Au} (p. 54-55, formula (13)). $\square$
\end{enumerate}
\medskip

We deduce from what preceeds an example of big actions with $G_2$ abelian of exponent $p^2$, with a small $p$-rank.
More precisely, we construct a subextension of $K_S^{m_2}$ with the commutative diagram:
$$
\begin{matrix}
0&\longrightarrow &G_S(m_2)&\longrightarrow& G(m_2)&\longrightarrow&\F_q&\longrightarrow &0\\
 &   & \varphi \, \downarrow&   &\downarrow&    &\vert \vert&    &  \\
  0&\longrightarrow &H&\longrightarrow& G&\longrightarrow&\F_q&\longrightarrow &0\\
    &   & \downarrow&   &\downarrow&    & &    &  \\
    &   & 0&   &0&    & &    &   
\end{matrix}   
$$
such that the pair $(C(m_2)/Ker(\varphi), G)$ is a big action where $G_2 \simeq \Z/p^2 \Z \times (\Z/p\Z)^t$ with $t=O(log_p \,g)$, $g$ being the genus of the curve $C(m_2)/Ker(\varphi)$.
Contrary to the previous case where the stability  under the translations by $\F_q$ was ensured by the maximality of $K_S^{m_2}$, the difficulty now lies in producing a system of equations defining a subextension of $K_S^{m_2}$  which remains globally invariant through the action of the group of translations $X \rightarrow X+y$, $y \in \F_q$. Write $q=p^e$. We have to distinguish the case $e$ even and $e$ odd.

\begin{proposition}
We keep the notation defined above. In particular, $K=\F_q(X)$ with $q=p^e$.
Assume $e=2\,s$, with $s \geq 1$, and put $r:=p^s$.
We define 
$$f_0(X):= a\, X^{1+r} \; \; \mbox{with} \; \; a \neq 0, \, \;  a \in \Gamma:=\{ \gamma \in \F_q, \gamma^r+\gamma=0\}$$
 and 
$$ \forall  \, i \in \{1,\cdots,p-1\}, \, \; f_i(X)=X^{ir/p} \, (X^q-X)=X^{ip^{s-1}}\, (X^q-X)$$ Let $L:=K(W_i)_{1 \leq i \leq p}$ be the extension of $K$ parametrized by the Artin-Schreier-(Witt) equations:
$$W_0^p-W_0=f_0(X) \quad \forall \, i\,  \in \{1,\cdots,p-1\}, \; W_i^q-W_i=f_i(X) \quad \mbox{and} \quad [W_0,W_p]^p-[W_0,W_p]= [f_0(X),0]$$
For all $i$ in $\{0,1,\cdots, p-1\}$, put $L_i:=K(W_0, \cdots,W_i)$.
\begin{enumerate}
\item $L$ is an abelian extension of $K$ such that every place in $S$ completely splits in $L$. Moreover,
$$L_0 \subset K_S^{r+2} \quad, \forall \, i \in \{1,\cdots,p-1\}, \,  \, L_i \subset K_S^{p^{s+1}+i+1} \,\, \mbox{with} \,  L \subset K_S^{m_2}$$
where $m_2=p^{s+1}+p+1$ is the integer defined in Proposition 6.4. (see table below).
\item The extension $L/K$ has degree $[L:K]=p^{2+(p-1)e}$.
Let $G_L$ be its Galois group. Then
$$G_L \simeq \Z/p^2\Z \times (\Z/p\Z)^t \quad \mbox{with} \, \, t=(p-1)\,e$$
\item The extension $L/K$ is stable under the translations: $X\rightarrow X+y$, with $y \in \F_q$. Therefore, the translations by $\F_q$ extend to form a $p$-group of $\F_q$-automorphisms of $L$, say $G$, with the following exact sequence:
 $$0\longrightarrow G_L\longrightarrow G \longrightarrow\F_q\longrightarrow 0$$
\item Let $g_L$ be the genus of the extension $L/K$. Then,
$$g_L= \frac{1}{2} \, \{ \, p^{2+2\,s\,(p-1)} \, (p^{s+1}+p-1)- p^{s} \, (p^2-p+1)-  p^{2\,s+1}\, (\sum_{i=0}^{p-2}\, q^i)\,  \}$$
In particular, when $e$ grows large, $g_L\sim  \frac{1}{2}\, p^{(2p-1)\, \frac{e}{2}+3}$ and $t=O(log_p \,g_L)$.\\
Note that, for $p=5$ and $e=4$, one gets $\frac{|G|}{g_L} \simeq 9,7049 \cdots$, which is slightly bigger than the quotient obtained for the whole extension $K_S^{m_2}$ (see \eqref{quotient}).
\end{enumerate}
\end{proposition}

\noindent \textbf{Proof:}
\begin{enumerate}
\item Fix $y$ in $\F_q$ and call $P_y:=(X-y)$, the corresponding place in $S$. As $f_i(y)=0$ for all $i$ in $\{1,\cdots, p-1\}$, the place $P_y$ completely splits in each extension $K(W_i)$ with $W_i^q-W_i=f_i(X)$. Therefore, to prove that $P_y$ completely splits in $L$, it is sufficient to show that $[f_0(y), 0] \in \wp(W_2(\F_q))$. By \cite{Bo} (Chap. IX, ex. 18),
this is equivalent to show that $Tr([f_0(y),0])=0$, where $Tr$ means the trace map from $W_2(\F_q)$ to $W_2(\F_p)$. We first notice that, when $y$ is in $\F_q$, $\gamma:=f_0(y)=a\,y^{1+r}$ lies in $\Gamma.$ It follows that:
$$Tr([\gamma,0])=\sum_{i=0}^{2s-1} \, F^i \, [\gamma,0]= \sum_{i=0}^{s-1} \, [\gamma^{p^i},0]+\sum_{i=0}^{s-1} \,[\gamma^{r\,p^i},0]=\sum_{i=0}^{s-1} \, [\gamma^{p^i},0]+\sum_{i=0}^{s-1} \,[-\gamma^{p^i},0]$$
For $p>2$, $[-\gamma^{p^i},0]=-[\gamma^{p^i},0]$ and $Tr([\gamma,0])=0$.
For $p=2$, since $p\,[\gamma,0]=[0,\gamma^p]$, one gets: 
$$Tr([\gamma,0])=[0,\gamma^p +\gamma^{p^2} +\cdots+\gamma^{p^s}]=[0,(\gamma +\gamma^{p} +\cdots+\gamma^{p^{s-1}})^p]= [0,Tr_{\F_r/\F_p} (\gamma)^p]$$
As $\Gamma$ coincides with $\{\beta^{r}-\beta, \, \beta \in \F_q\}$ (see e.g. \cite{Au} p. 58), $Tr_{\F_r/\F_p} (\gamma)=0$ and $Tr([\gamma,0])=0$.
To establish the expected inclusions, it remains to compute the conductor of each extension $L_i$. First of all, \cite{Au} (I, ex. 3.3) together with \cite{St93} (Prop III,7.10) shows that the conductor of $L_0$ is $r+2$. Thus, $L_0 \subset K_S^{r+2}$. Moreover, as $f_i(X)=X^{i+p^{s+1}}-X^{1+ip^{s-1}}$ mod $\wp(\F_q[X])$, we infer from \cite{Au} (I, ex. 3.3) and \cite{Au} (I, Cor. 3.2) that the conductor of $L_i$ is $1+i+p^{s+1}$. So, $L_i \subset K_S^{1+i+p^{s+1}}$. To complete the proof, it remains to show that $L$ has conductor $m_2$, which derives from \cite{Ga} (see comments above).
\medskip

The equations, conductor and degree of each extension $L_i$ are finally gathered in the table below.
\medskip

\begin{tabular}{|c|c|c|c|}
\hline
$L_i$ & conductor $m$ & $[L_i:K]$ & New equations\\
\hline
$K$ & $0 \leq m \leq r+1=p^s+1$ & $1$ & \\
\hline
$L_0$ & $r+2 \leq m \leq p^{s+1}+1=m_2-p$ & $p$ & $W_{0}^p-W_{0}=f_0(X)$\\
\hline
$L_1$ & $m=p^{s+1}+2=m_2-(p-1)$ & $p^{1+e}$ & $W_1^q-W_1=f_1(X)$\\
\hline
$L_2$ & $m=p^{s+1}+3=m_2-(p-2)$ & $p^{1+2e}$ & $W_2^q-W_2=f_2(X)$\\
\hline
$\cdots \cdots$ & $\cdots \cdots$ & $\cdots \cdots$ & $\cdots \cdots$\\
\hline
$L_i$  & $m=p^{s+1}+i+1=m_2-(p-i)$ & $p^{1+ie}$ & $W_i^q-W_i=f_i(X)$\\
\hline
$\cdots \cdots$ & $\cdots \cdots$ & $\cdots \cdots$ & $\cdots \cdots$\\ 
\hline
$L_{p-1}$ & $m=p^{s+1}+p=m_2-1$ & $p^{1+(p-1)e}$ & $W_{p-1}^q-W_{p-1}=f_{p-1}(X)$\\
\hline
$L$ & $m=p^{s+1}+p+1=m_2$ & $p^{2+(p-1)e}$ & $[W_0,W_p]^p-[W_0,W_p]=[f_0(X),0] $\\
\hline
\end{tabular}

\medskip

\item See table above.

\item Fix $y$ in $\F_q$. Consider $\sigma$ in $G(m_2)$ (defined as in Proposition 6.5) such that $\sigma(X)=X+y$. \begin{enumerate}
\item We first prove that $\sigma(W_0) \in L_0$. Indeed, as $y \in \F_q$ and $a \in \Gamma=\{ \gamma \in \F_q, \gamma^r+\gamma=0\}$, 
$$ \begin{array}{ll}
\wp(\sigma(W_0)-W_0)&= \sigma(\wp(W_0))-\wp(W_0)\\
& =f_0(X+y)-f_0(X)\\
&= a\,y\,X^r+a\,y^r\,X+f_0(y)\\
&=-a^r \,y^{r^2} \, X^r+a\,y^r\,X+f_0(y)\\
&=\wp(P_y(X))+f_0(y)
\end{array}$$
where $P_y(X):=(I+F+F^2+\cdots +F^{s-1}) \, (-a\, y^r \, X)$. Since $f_0(y) \in \wp(\F_q)$ (see proof of the first point), it follows that $\wp(P_y(X))+f_0(y)$ belongs to $\wp(\F_q[X])$.
Therefore, $\sigma(W_0) \in L_0=\F_q(X,W_0)$.
\item We now prove that, for all $i$ in $\{1,\cdots,p-1\}$, $\sigma(W_i) \in L_i$. Indeed,
$$
\begin{array}{ll}
(F^e-\id)\, (\sigma(W_i)-W_i)&= \sigma (W_i^q-W_i)-(W_i^q-W_i)\\
\quad \\
&=f_i(X+y)-f_i(X)\\
\quad \\
&=(X+y)^{i \,p^{s-1}} \, (X^q-X)-X^{i\, p^{s-1}}\, (X^q-X)\\
\quad \\
&=(X^{p^{s-1}}+y^{p^{s-1}})^{i} \, (X^q-X)-X^{i\, p^{s-1}}\, (X^q-X)\\
\quad \\
&=\sum_{j=1}^{i-1} \, \binom ij \, y^{(i-j)p^{s-i}}\, f_j(X) \, \mod (F^e-\id) \, (\F_q[X])\\ 
\quad \\
&=(F^e-\id) \, (\sum_{j=1}^{i-1} \, \binom ij \, y^{(i-j)p^{s-i}}\, W_j) \, \mod (F^e-\id) \, (\F_q[X])\\ 
\end{array}$$
where the sum is empty for $i=1$. It follows that $\sigma(W_i) \in L_i=\F_q(X,W_0, W_1, \cdots, W_i)$.
\item To conclude, we show that $\sigma(W_p) \in L$, which requires the use of Remark 6.3.4. Indeed, compute:
$$ \begin{array}{ll}
\Delta&:=\wp(\sigma \, [W_0,W_p]-[W_0,W_p])\\
&= \sigma (\wp([W_0,W_p])-\wp([W_0,W_p])\\
&=[f_0(X+y),0]-[f_0(X),0]
\end{array}$$
As shown in the proof of the first point, $[f_0(y),0]$ lies in $\wp(W_2(\F_q))$. Then, 
$$\Delta=[f_0(X+y),0]-[f_0(X),0]-[f_0(y),0]-[P_y(X),0]+[P_y(X),0]^p \mod \,\wp(W_2(\F_q[X]))$$
with $y$ in $\F_q$ and $P_y$ defined as above. 
Let $W(\F_q)$ be the ring of Witt vectors with coefficients in $\F_q$. Then, for any $y \in \F_q$, we denote by $\tilde{y}$ the Witt vector
$\tilde{y}:=(y,0,0,\cdots) \in W(k)$. For any $P(X):=\sum_{i=0}^{s} \, a_i \, X^i \in \F_q[X]$, we denote by $\tilde{P}(X):=\sum_{i=0}^s \, \tilde{a_i} \, X^i \in W(\F_q)[X]$. 
The addition in the ring of Witt vectors yields: 
$$ \begin{array}{ll}
\Delta&= [0,A]  \quad \mod \,\wp(W_2(\F_q[X]))
\end{array}$$
 where  $A$ is the reduction modulo $p\, W_2(\F_q)[X]$ of:
$$\frac{1}{p} \{\tilde{f_0}(X+\tilde{y})^p-\tilde{f_0}(X)^p-\tilde{f_0}(\tilde{y})^p+\tilde{P}_y(X)^p-\tilde{P}_y(X)^{p^2} -(\tilde{f_0}(X+\tilde{y})-\tilde{f_0}(X)-\tilde{f_0}(\tilde{y})-\tilde{P}_y(X)+\tilde{P}_y(X)^p)^p\} $$
Since $\tilde{f_0}(X+\tilde{y})-\tilde{f_0}(X)-\tilde{f_0}(\tilde{y})+\tilde{P}_y(X)-\tilde{P}_y(X)^p=0$ mod $p\, W(\F_q)[X]$, $A$ becomes:
$$A=\frac{1}{p} \, \{\tilde{f_0}(X+\tilde{y})^p-\tilde{f_0}(X)^p-\tilde{f_0}(\tilde{y})^p+\tilde{P}_y(X)^p-\tilde{P}_y(X)^{p^2}\} \quad \mod \, p\,W(\F_q)[X]  $$
We observe that:
$$
\begin{array}{lll}
\tilde{f_0}(X+\tilde{y})^p&=\tilde{a}^p \, (X+\tilde{y})^p \, (X+\tilde{y})^{p^{s+1}}& \mod \, p^2\,W(\F_q)[X]\\
\quad \\
&= \tilde{a}^p \, (X+\tilde{y})^p \, (X^{p^s}+\tilde{y}^{p^s})^{p} \,& \mod \, p^2\,W(\F_q)[X]\\
\quad \\
&=\tilde{a}^p \, \sum_{i=0}^p \,\sum_{j=0}^p \, \binom pi \, \binom pj \, X^{j+ip^s} \, \tilde{y}^{p-j+p^s\,(p-i)}  \, &\mod \, p^2\,W(\F_q)[X]\\
\end{array} $$
As $\binom pi \, \binom pj=0$ mod $p^2$ when $0<i<p$ and $0<j<p$, one obtains:
$$
\tilde{f_0}(X+\tilde{y})^p-\tilde{f_0}(X)^p-\tilde{f_0}(\tilde{y})^p=\tilde{a}^p \, \sum_{(i,j) \in I}
\, \binom pi \, \binom pj \, X^{j+ip^s} \, \tilde{y}^{p-j+p^s\,(p-i)}  \quad \mod \, p^2\,W(\F_q)[X]\\
 $$
with 
$$I:=\{(i,j) \in \N^2, \, 0 \leq i \leq p, 0\leq j \leq p,\,  ij=0 \mod p, (i,j)\neq(0,0) \,, (i,j)\neq(p,p)\}$$
Besides, 
$$\begin{array}{lll}
\tilde{P}_y(X)^p-\tilde{P}_y(X)^{p^2}&= (\sum_{i=0}^{s-1} \, (-\tilde{a} \, \tilde{y}^r \, X)^{p^i}) ^{p} -
(\sum_{i=0}^{s-1} \, (-\tilde{a} \, \tilde{y}^r \, X)^{p^i}) ^{p^2} &\, \mod \, p^2\,W(\F_q)[X] \\
\quad \\
&= (\sum_{i=0}^{s-1} \, (-\tilde{a} \, \tilde{y}^r \, X)^{p^i})^{p} -
(\sum_{i=0}^{s-1} \, (-\tilde{a} \, \tilde{y}^r \, X)^{p^{i+1}})^{p} &\, \mod \, p^2\,W(\F_q)[X] \\
\quad \\
&= -\tilde{a}^p\, \tilde{y}^{rp} \, X^p+\tilde{a}^{rp} \, \tilde{y}^{r^2p} \, X^{pr} +p \, \tilde{T}_y(X) & \mod \, p^2\,W(\F_q)[X] \\
\end{array}$$
with $\tilde{T}_y(X) \in W(\F_q)[X]$. As $y \in \F_q$ and $a\in \Gamma$, we get:
$$\begin{array}{lll}
\tilde{P}_y(X)^p-\tilde{P}_y(X)^{p^2}&= -\tilde{a}^p\, \tilde{y}^{rp} \, X^p-\tilde{a}^{p} \, \tilde{y}^{p} \, X^{pr} +p \, \tilde{T}_y(X) & \mod \, p^2 \,W(\F_q)[X]\\
\end{array}$$
As a consequence, 
$$A=\tilde{a}^p \, \sum_{(i,j)\in I_1}
\, \frac{1}{p} \, \binom pi \, \binom pj \, X^{j+ip^s} \, \tilde{y}^{p-j+p^s\,(p-i)} +\tilde{T}_y(X) \quad \mod \, p\, \wp(\F_q[X])$$
with $$I_1:=\{(i,j) \in I,\, (i,j)\neq(0,p) \,, (i,j)\neq(p,0)\}$$ So, $A$ reads:
$$A=a^p \, \sum_{(i,j)\in I_1}
\, \frac{1}{p} \, \binom pi \, \binom pj \, X^{j+ip^s} \, y^{p-j+p^s\,(p-i)} +T_y(X) $$
with $T_y \in \F_q[X]$. We first consider the sum. 
Since, for $i=0$, $i=p$ and $j=p$, one gets monomials whose degree (after eventual reduction mod $\wp(\F_q[X])$) is strictly lower than $1+p^{s}$, one can write:
$$A=a^p \,\sum_{j=1}^{p-1} \, \frac{1}{p} \, \binom pj \, X^{j+ip^s} \, y^{p-j} +R_y(X)+T_y(X)  \quad\,  \mod \, \wp(\F_q[X])$$
where $R_y(X)$ is a polynomial of $\F_q[X]$ with degree strictly lower than $1+p^s=1+r$.
We now focus on the polynomial $T_y(X) \in \F_q[X]$. It is made of monomials which read either 
$X^{i_0+i_1\,p+\cdots +i_{s-1}\,p^{s-1}}$ with $i_0+i_1+\cdots+i_{s-1}=p$ or $X^{i_1\,p+\cdots+i_{s} \, p^s}$,
 with $i_1+i_2+\cdots+i_{s}=p$. As $X^{i_1\,p+\cdots+i_{s} \, p^s}=X^{i_1+\cdots+i_{s} \, p^{s-1}}$ mod $\wp(\F_q[X])$, it follows that $T_y$ does not have any monomial with degree higher than $1+p^s$ after reduction mod $\wp(\F_q[X])$. So,
$$A=a^p \,\sum_{j=1}^{p-1} \,\frac{1}{p} \, \binom pj \, X^{j+ip^s} \, y^{p-j} +R^{[1]}_y(X)  \quad\,  \mod \, \wp(\F_q[X])$$
where $R^{[1]}_y(X)$ is a polynomial of $\F_q[X]$ with degree strictly lower than $1+r$.
Since, for all $j$ in $\{1,\cdots,p-1\}$,  $f_j(X)= X^{j+p^{s+1}}-X^{1+jp^{s-1}}$ mod $\wp(\F_q[X])$, we gather: 
$$A=a^p \,\sum_{j=1}^{p-1} \,\frac{1}{p} \, \binom pj  \, y^{p-j}\, f_j(X) +R^{[2]}_y(X)  \quad\,  \mod \, \wp(\F_q[X])$$
where $R^{[2]}_y(X)$ is a polynomial of $\F_q[X]$ with degree strictly lower than $1+r$. Then,
$$
\begin{array}{lll}
A &= \sum_{j=1}^{p-1} \, c_j(y) \, f_j(X)+R^{[2]}_y(X)  \, & \mod \, \wp(\F_q[X])\, 
\end{array}$$
with $c_j(y):=a^p\, \frac{1}{p} \, \binom pj  \, y^{p-j} \in \F_q$.
It follows that:
$$
\begin{array}{lll}
A&= \sum_{j=1}^{p-1} \, (F^e-\id) \, (c_j(y) \, W_j) +R^{[2]}_y(X)  \,  &\mod \, \wp(\F_q[X]) \\
&\\
&= (F-\id)\, \sum_{j=1}^{p-1} \,  \, P_j( W_j) +R^{[2]}_y(X)  \, & \mod \, \wp(\F_q[X]) 
\end{array}$$
where $P_j(W_j)= (\id +F+\cdots +F^{e-1})\, (c_j(y)\, W_j) \in \F_q[W_j]$. 
We gather that:
$$\wp(\sigma \, [W_0,W_p]-[W_0,W_p])=\wp \, ([0,\sum_{j=1}^{p-1} \,  \, P_j( W_j)])+ [0, R^{[2]}_y(X)] \, \mod \,\wp(W_2(\F_q[X])) $$
As a consequence, $[0, R^{[2]}_y(X)]$ lies in $\wp(W_2(K_S^{m_2}))$ and so, there exists $V \in K_S^{m_2}$ such that $V^p-V=R^{[2]}_y(X)$
Accordingly, $K(V)$ is a $K$-subextension of $K_S^{m_2}$ with conductor $1+deg(R^{[2]}_y(X)) \leq 1+r$. In particular, $K(V) \subset K_S^{r+1}=K=\F_q(X)$, which implies that $R^{[2]}_y(X)\in \wp(K)$.
Therefore, $$\wp(\sigma \, [W_0,W_p]-[W_0,W_p])=\wp  \, ([0,\sum_{j=1}^{p-1}   \, P_j( W_j)]) \, \mod \,\wp(W_2(K))$$ which allows to conclude that $\sigma \, (W_p)$ is in $L=K(W_0,W_1,\cdots,W_p)$.
\end{enumerate}

\item As $L \subset K_S^{m_2}$ and $L \not \subset K_S^{m_2-1}$, the formula established in Proposition 6.5.3. yields:
$$
\begin{array}{ll}
g_L&=1+[L:K] \,(-1+\frac{m_2}{2} )-\frac{1}{2} \sum_{j=0}^{m_2-1} [K_S^{j}:K] \\
&\\
&=1+p^{2+(p-1)e}\, (-1+\frac{p^{s+1}+p+1}{2})-\frac{1}{2}(r+2+(m_2-p-(r+2)+1)\, p+
\sum_{i=1}^{p-1} \, p^{1+i\,e} )\\
&\\
&= \frac{1}{2} \, p^{2+(p-1)e} \, (p^{s+1}+p-1) -\frac{1}{2} \, (p^s+p^{s+2}-p^{s+1}+
\sum_{i=1}^{p-1} \, p^{1+i\,2\,s}) \\
&\\
&=\frac{1}{2} \, p^{2+(p-1)e} \, (p^{s+1}+p-1) -\frac{1}{2} \,p^s (p^2-p+1) -\frac{1}{2} \, p^{2s+1} (1+q+q^2+\cdots +q^{p-2}) \qquad \square
\end{array}
$$
\end{enumerate}

The preceding proposition can be generalized to construct a big action endowed with a second ramification group $G_2$ abelian of exponent as large as we want.

\begin{proposition}
We keep the notation defined above. In particular, $q=p^e$, with $e=2s$ and $s\geq 1$.
Let $n \geq 2$. Put $m_n:=1+p^{n-1}\, (1+p^s)$.
If $\frac{q}{-1+m_n/2} >\frac{2\,p}{p-1}$, the pair $(C(m_n), G(m_n))$, as defined in Proposition 6.5, is a big action with a second ramification group $G_S(m_n)$ abelian of exponent at least $p^n$.
\end{proposition}

\noindent \textbf{Proof:} Proposition 6.5.4 first ensures that $(C(m_n), G(m_n)$ is a big action.
Consider the $p^n$-cyclic extension $K(W_1,\cdots, W_n)/K$ parametrized as follows, with Witt vectors of length $n$:
$$[W_1,\cdots,W_n]^p-[W_1,\cdots,W_n]=[f_0(X), 0 ,\cdots, 0]$$
where $f_0(X)=a\,X^{1+r}$ is defined as in Proposition 6.7, i.e. $r=p^s$, $a^r+a=0$ , $a\neq 0$. The same proof as in Proposition 6.7.1 shows that all places of $S$ completely split in $K(W_1,\cdots,W_n)$. Moreover, by \cite{Ga} (Thm. 1.1) the conductor of the extension $K(W_1,\cdots,W_n)$ is $1+max \{p^{n-1}\, (1+p^s),0\}=m_n$. It follows that $K(W_1,\cdots, W_n)$ is included in $K_S^{m_n}$. Therefore, $G_S(m_n)$ has a quotient of exponent $p^n$ and the claim follows. $\square$

\medskip

Tne next proposition is an analogue of Proposition 6.7 in the case where $e$ is odd.
We does not mention the proof which is mainly similar to the proof of Proposition 6.7.

\begin{proposition}
We keep the notation defined above. In particular, $K=\F_q(X)$ with $q=p^e$.
Assume $e=2\,s-1$, with $s \geq 2$, and put $r:=\sqrt{qp}=p^{s}$.
We define 
$$ \forall  \, i \in \{1,\cdots,p-1\}, \; f_i(X)=X^{ir/p} \, (X^q-X)=X^{ip^{s-1}}\, (X^q-X)$$ $$\forall  \, i \in \{1,\cdots,p-1\}, \; g_i(X)=X^{ir/{p^2}} \, (X^q-X)=X^{ip^{s-2}}\, (X^q-X)$$ 
Let $L:=K(W_i,V_j)_{1\leq i \leq p, 1 \leq j \leq p-1 }$ be the extension of $K$ parametrized by the Artin-Schreier-(Witt) equations:
$$\forall \, i \in \{1,\cdots,p-1\}, \; W_i^q-W_i=f_i(X) \quad \mbox{and} \quad
\forall \, j\,  \in \{1,\cdots,p-1\}, \; V_j^q-V_j=g_j(X)$$
$$[W_1,W_p]^p-[W_1,W_p]= [X^{1+p^s},0]-[X^{1+p^{s-1}},0]$$
For all $i$ and $j$ in $\{1,\cdots, p-1\}$, put $L_{i,0}:=K(W_k)_{1\leq k \leq i}$
and $L_{p-1,j}:=K(W_i,V_k)_{1\leq i \leq p-1, 1 \leq k \leq j}$.
\begin{enumerate}
\item $L$ is an abelian extension of $K$ such that every place in $S$ completely splits in $L$.
Then, 
$$\forall \, i, \, j   \in \{1,\cdots,p-1\}, \, L_{i,0} \subset K_S^{p^s+i+1} \quad, \, \, L_{p-1,j} \subset K_S^{p^{s+1}+j+1} \quad  \mbox{and}  \quad  L \subset K_S^{m_2}$$
where $m_2=p^{s+1}+p+1$ is the integer defined in Proposition 6.4. (see table below.)

\item The extension $L/K$ has degree $[L:K]=p^{2(p-1)e+1}$.
Let $G_L$ be its Galois group. Then
$$G_L \simeq \Z/p^2\Z \times (\Z/p\Z)^t \quad  \mbox{with} \,\; t=2\,(p-1)\,e-1$$
\item The extension $L/K$ is stable under the translations: $X\rightarrow X+y$, with $y \in \F_q$. Therefore, the translations by $\F_q$ extend to form a $p$-group of $\F_q$-automorphisms of $L$, say $G$, with the following exact sequence:
 $$0\longrightarrow G_L\longrightarrow G \longrightarrow\F_q\longrightarrow 0$$
\item Let $g_L$ be the genus of the extension $L/K$. Then,
$$g_L= \frac{1}{2} \, \{ \, p^{1+(2p-1)e} \, (p^{s+1}+p-1)- p^{(p-1)e} \, (p^{s+1}-p^s-p+1)- p^s+p^e\, (\sum_{i=0}^{2p-3} q^i) \}$$
In particular, when $e$ grows large, $g_L\sim  \frac{1}{2}\, p^{2+4s(p-1)+s}$ and $t=O(log_p \, g_L)$.
\end{enumerate}
\end{proposition}

We gather in the table below the conductors, degrees and equations of each extension.
\medskip

\begin{tabular}{|c|c|c|c|}
\hline
$L_{i,j}$ & conductor $m$  &$[L_{i,j}:K]$ &  New equations\\
\hline
$K$& $0 \leq m \leq r+1=p^s+1$ &  $1$ & \\
\hline
$L_{1,0}$ & $m=r+2=p^s+2 $ & $p^e$ & $W_{1}^q-W_{1}=f_1(X)$\\
\hline
$\cdots \cdots$ & $\cdots \cdots$& $\cdots \cdots$ & $\cdots \cdots$ \\
\hline
$L_{i,0} $ & $m=p^{s}+i+1$ & $p^{ie}$ &  $W_i^q-W_i=f_i(X)$\\
\hline
$\cdots \cdots$ & $\cdots \cdots$ & $\cdots \cdots$& $\cdots \cdots$\\ 
\hline
$L_{p-1,0}$& $p^{s}+p \leq m \leq p^{s+1}+1$  & $p^{(p-1)e}$ &  $W_{p-1}^q-W_{p-1}=f_{p-1}(X)$\\
\hline
$L_{p-1,1}$ & $m=p^{s+1}+2=m_2-(p-1) $ & $p^{pe}$ & $V_{1}^q-V_{1}=g_1(X)$\\
\hline
$\cdots \cdots$ & $\cdots \cdots$& $\cdots \cdots$ & $\cdots \cdots$ \\
\hline
$L_{p-1,j}$ & $m=p^{s+1}+j+1=m_2-(p-j)$ & $p^{(p+j-1)e}$ &  $V_j^q-V_j=g_j(X)$\\
\hline
$\cdots \cdots$ & $\cdots \cdots$ & $\cdots \cdots$ & $\cdots \cdots$\\ 
\hline
$L_{p-1,p-1}$ & $ m=p^{s+1}+p=m_2-1$  & $p^{2(p-1)e}$ &  $V_{p-1}^q-V_{p-1}=g_{p-1}(X)$\\
\hline
$L$ & $m=p^{s+1}+p+1=m_2$ & $p^{1+2\,(p-1)e}$ &  $[W_1,W_p]^p-[W_1,W_p]=$\\
& & & $[X^{1+p^s},0]-[X^{1+p^{s-1}},0] $\\
\hline
\end{tabular}
\medskip

\section{A local approach to big actions.}

\indent Let $(C,G)$ be a big action. We recall that there exists a point $\infty \in C$ such that $G$ is equal to $G_1(\infty)$ the wild inertia subgroup of $G$ at $\infty$, which means that the cover $\pi: C \rightarrow C/G$ is totally ramified at $\infty$. 
Moreover, the quotient curve $C/G$ is isomorphic to the projective line: $\p_k^1$ and $\pi$ is \'etale above the affine line: $\A_k^1=\p_k^1-\pi(\infty)=Spec\, k[T]$. The inclusion $k[T] \subset k((T^{-1}))$ induces a Galois extension $k(C) \otimes_{k(T)} k((T^{-1}))=:k((Z))$ over $k((T^{-1}))$
with group equal to $G$ and ramification groups in lower notation equal to $G_i:=G_i({\infty})$.
Then, the genus of $C$ is given by the formula: $g=\frac{1}{2}\, (\sum_{i\geq 2} (|G_i|-1)) >0$ (see \eqref{genus}).
It follows that: $$\frac{|G|}{\sum_{i\geq 2} (|G_i|-1)}=\frac{|G|}{2\, g} > \frac{p}{p-1}.$$ This leads to the definition below.

\begin{defi}
We call "local big action" any pair $(k((Z)),G)$ where  $G$ is a finite $p$-subgroup of
$Aut_k(k((Z))$ whose ramification groups in lower notation at $\infty$ satisfy the two inequalities:  $$g(G):=\frac{1}{2} (\sum_{i\geq 2} (|G_i|-1))>0 \qquad \mbox{and} \qquad \frac{|G|}{g(G)} > \frac{2\: p}{p-1}.$$
\end{defi}

It follows from the Katz-Gabber Theorem (see \cite{Ka} Thm. 1.4.1 or \cite{PGi} cor. 1.9) that big actions $(C,G)$ and local big actions $(k((Z)),G)$ are in 1-to-1 correspondance via the following functor induced by the inclusion $k[T] \subset k((T^{-1}))$:

$$\left\{
\aligned \mbox{finite \'etale  Galois covers of Spec  k[T]} \\
\mbox{ with Galois group a p-group}
 \endaligned
\right\} \quad \longrightarrow \quad \left\{
\aligned \mbox{finite \'etale  Galois covers of  Spec } k((T^{-1}))\\
 \mbox{with  Galois  group a  p-group} \endaligned
\right\}
$$

Therefore, we can infer from the global point of view properties related to local extensions that would be difficult to prove directly.
For instance, if $(k((Z)), G)$ is a local big action, we can deduce that $G_2$ is stricly included in $G_1$. Furthermore, we obtain:
$$\frac{|G|}{g(G)^2} \leq \frac{4\, p}{(p-1)^2}.$$

%%\addcontentsline{toc}{chapter}{Bibliographie}

\medskip

\begin{flushleft}

Michel MATIGNON \\
Institut de Math\'ematiques de Bordeaux, 
Universit\'e de Bordeaux I,
351 cours de la Lib\'eration, 
33405 Talence Cedex, France \\
e-mail : {\tt Michel.Matignon@math.u-bordeaux1.fr}

\bigskip

Magali ROCHER\\
Institut de Math\'ematiques de Bordeaux,
Universit\'e de Bordeaux I,
351 cours de la Lib\'eration, 
33405 Talence Cedex, France \\
e-mail : {\tt Magali.Rocher@math.u-bordeaux1.fr}

\end{flushleft}

\end{document}